\documentclass[a4pape]{cas-dc}
\usepackage[numbers]{natbib}

\usepackage{color}
\graphicspath{{img/}}
\usepackage{amsfonts}
\usepackage{graphicx,psfrag}
\usepackage{amsmath,bm,amssymb}
\usepackage{hyperref}
\usepackage{booktabs}
\usepackage{multirow}  
\usepackage{makecell} 
\allowdisplaybreaks[4]

\usepackage{multicol}

\def\tsc#1{\csdef{#1}{\textsc{\lowercase{#1}}\xspace}}
\tsc{WGM}
\tsc{QE}
\tsc{EP}
\tsc{PMS}
\tsc{BEC}
\tsc{DE}
\newtheorem{definition}{Definition}
\newtheorem{assump}{\bf{Assumption}}
\newtheorem{theorem}{Theorem}
\newtheorem{lemma}{Lemma}
\newdefinition{rmk}{Remark}
\newproof{proof}{Proof}

\newproof{pot}{Proof of Theorem \ref{thm2}}

\newdefinition{prop}{Proposition}

\begin{document}
\let\WriteBookmarks\relax
\def\floatpagepagefraction{1}
\def\textpagefraction{.001}
\shorttitle{ }
\shortauthors{J. Qi, J. Zhang, M. Krstic}
	
\title [mode = title]{Neural Operators for PDE Backstepping Control of First-Order Hyperbolic PIDE with Recycle and Delay}
\tnotemark[1]
	
\tnotetext[1]{ This work was partially supported by the National Natural Science Foundation of China (62173084), the Project of Science and Technology Commission of Shanghai Municipality, China (23ZR1401800, 22JC1401403).}
	

\author[DHU,Textile]{Jie Qi}
\cormark[1]
\ead{jieqi@dhu.edu.cn}
	
\author[DHU]{Jing Zhang} 
\ead{zhangjing@mail.dhu.edu.cn}
	
\author[UCSD]{Miroslav Krstic}
\ead{krstic@ucsd.edu}
	
\address[DHU]{College of Information Science and Technology, Donghua University, Shanghai, 201620, China}
\address[Textile]{Engineering Research Center of Digitized Textile and Fashion Technology Ministry of Education, Donghua University, Shanghai, 201620, China} 
\address[UCSD]{Dept. of Mechanical and Aerospace Engineering, University of California San Diego, La Jolla, CA 92093-0411, USA.}

\cortext[cor1]{Corresponding author}
	%
	
	\begin{abstract}The recently introduced DeepONet operator-learning framework for PDE control is extended from the  results for basic hyperbolic and parabolic PDEs to an advanced hyperbolic class that involves delays on both the state and the system output or input. The PDE backstepping design produces gain functions that are outputs of a nonlinear operator, mapping functions on a spatial domain into functions on a spatial domain, and where this gain-generating operator's inputs are the PDE's coefficients. The operator is approximated with a DeepONet neural network to a degree of accuracy that is provably arbitrarily tight. Once we produce this  approximation-theoretic result in infinite dimension, with it we establish stability in closed loop under feedback that employs approximate gains. 
In addition to supplying such results under full-state feedback, we also develop DeepONet-approximated observers and output-feedback laws and prove their own stabilizing properties under neural operator approximations. With numerical simulations we illustrate the theoretical results and quantify the numerical effort savings, which are of two orders of magnitude, thanks to replacing the numerical PDE solving with the DeepONet. 

	\end{abstract}
	
	
	
 \ExplSyntaxOn
 \keys_set:nn { stm / mktitle } { nologo }
 \ExplSyntaxOff
\maketitle
	
\section{Introduction}
\label{1}
	
	In \cite{bhan2023neural,krstic2023neural}, a method was introduced to pre-train the backstepping methodology, offline and once and for all, for certain entire classes of PDEs so that the implementation of the controller to any specific PDE within the class is nothing more than a function evaluation of a neural network  that produces the controller gains based on the specific plant coefficients of the PDE being controlled. 
	
	In this paper, we extend this method to a broader and more advanced class of hyperbolic partial integro-differential systems, which involve delays on the state and the output or input.
	
	\subsection{The broader context of learning-based and data-driven control}
	
	Recently, learning-based control approaches have  attracted great attention due to their leveraging of capabilities of deep neural networks. Some of these approaches learn control strategies from data without explicit knowledge of system dynamics, and some are able to deal with uncertainties and disturbances. Stability and robustness can be proven with some of these control methods \cite{jiang2020learning}, which builds trust for their use in practice.  Progress has taken place with learning-based model predictive control (MPC) for uncertain models \cite{soloperto2018learning, nguyen2021robust}, Lyapunov functional based control design \cite{abate2020formal, zhou2022neural}, reinforcement learning  (RL) based linear quadratic regulator \cite{mohammadi2021convergence, hambly2021policy}, and other methods.  
	RL has also been  applied to PID  tuning \cite{dogru2022reinforcement, mcclement2022meta}, with a notable use in \cite{lawrence2022deep}, where a deep RL-based PID tuning method is proposed and experimented on the physical two-tank system without prior pre-training.
	For the risks that might arise during the RL based control process, recently, safe reinforcement learning  has emerged as a new research focus, see e.g., \cite{garcia2015comprehensive, qin2021density, zhang2020first}.
	
	Learning-based control in  unmanned systems is pursued
	in \cite{chen2019knowledge, zheng2023ddpg, li2021maneuvering, wang2021data}. 
	For example in control of the unmanned aerial vehicles (UAVs), a dual-stream Actor-Critic network structure is applied to extract environmental features, enabling UAVs to safely navigate in environments with multiple obstacles\cite{zhang2022autonomous}. Data-driven control methods extract the hidden patterns from a large amount of data, which improves control performance in uncertain environment. In \cite{o2022neural}, a deep network learning-based trajectory tracking controller, called Neural-Fly, is proposed for drones' agile flight in rapidly changing strong winds. 
	Transfer learning also used to leverage control strategies and models that have already been learned to accelerate the learning and adaptation process for new tasks, e.g., \cite{venturini2021distributed, chu2022joint, li2021maneuvering}.

	
	\subsection{Learning-enhanced PDE control}
	
	Many engineering problems are spatio-temporal processes, often modeled by partial differential equations (PDEs) instead of ordinary differential equations (ODEs), such as plug flow reactor \cite{xu2018optimal}, traffic flow \cite{qi2022delay, yu2022simultaneous}, hydraulics and river dynamics \cite{bastin2019boundary}, pipeline networks \cite{aamo2015leak,anfinsen2022leak}, melt spinning processes \cite{gotz2010optimal}, flexible robots \cite{he2018pde}, flexible satellite \cite{he2015dynamic}, tokamaks \cite{mavkov2017distributed} and so on.  

PDE backstepping has been particularly effective in the stabilization of PDEs. Since this paper is focused on a hyperbolic  partial integro differential equation (PIDE)  class, we mention only a few designs for hyperbolic systems here. A design for a single hyperbolic PIDE was introduced in \cite{krstic2008}. A pair of coupled hyperbolic PDEs was stabilized next, with a single boundary input in \cite{Coron2013Local}. An extension to $n+1$ hyperbolic PDEs with a single input was introduced in \cite{di2013stabilization},  an extension to cascades with ODEs in \cite{DIMEGLIO2018281}, an extension to “sandwiched” ODE-PDE-ODE systems in \cite{WANG2020109131, 9319184}, and redesigns robust to delays in \cite{Auriol2018Delay,Auriol2018Delay1}.

Since the dynamics of the PDE systems are defined in infinite-dimensional function spaces, the gains in the PDE control systems (feedback controllers, observers, identifiers) are not vectors or matrices but functions of spatial arguments.
	When the coefficients of the system are spatially-varying, the equations governing the control gain kernels usually cannot be solved explicitly, as they are complex PDEs  
	and need to be solved numerically, e.g. \cite{di2013stabilization, hu2015control,vazquez2016boundary,vazquez2023kernel}. 
	When any coefficient changes, the control gain PDEs need to be re-solved, which is burdensome even if performed offline and once, let alone if it needs to be performed repeatedly in real time, in the context of adaptive control or gain scheduling. 
	
Operator learning  refers to the learning of an infinite-dimensional mapping operator by means of deep neural network.    It is of interest to find a neural network (NN) which learns control gain operators from a large set of previously offline-solved control design problems for a sample set of PDEs in a certain class. For example, \cite{margenberg2023optimal} utilized the Fourier Neural Operator to address the optimal Dirichlet boundary control problem in nonlinear optics. A robust framework employing an operator learning technique for such problems with PDEs constraints is provided in \cite{hwang2022solving}. Furthermore, \cite{lin2021operator}  demonstrates the application of DeepONet in learning the relationship between liquid pressure and bubble generation, thereby validating the efficacy and precision of neural operators (NO) in predicting the dynamics of multi-rate bubble growth.\color{black}
 
 The DeepONet framework \cite{lu2022multifidelity, bhan2023neural, krstic2023neural} is an efficient method for PDE control, because it not only speeds up computation, e.g., on the order of magnitude of $10^3$ time \cite{krstic2023neural}, as compared to solving for the control gains numerically, but also provides a methodology for stability analysis. The DeepONet \cite{lu2019deeponet} consists of two sub-networks, i.e.,  branch net and trunk net. The branch net encodes the discrete input function space and the trunk encodes the domain of the output functions. The combination of branch and trunk nets improves the generalization and efficiency of operators learning of the DeepONet, so that it realize the regression of infinite-dimensional functions from a relatively small number of datasets \cite{deng2022approximation}, which brings new insights for the learning based control. 
Furthermore, the universal approximation theorem \cite{lu2019deeponet, deng2022approximation, chen1995universal}, which states that a nonlinear continuous operator can be approximated by an appropriate DeepONet with any given approximation error, provides the basis for the rigorous stability analysis of the closed-loop system under the neural operator controller.

	This offline learning PDE control design framework was pioneered in \cite{bhan2023neural}. Among the PDE control design approaches, PDE backstepping was used, due to its non-reliance on model reduction and its avoidance of numerically daunting operator Riccati equations. Among the neural operator methods, the DeepONet \cite{lu2019deeponet, deng2022approximation} approach was employed, due to its availability of universal approximation theorem in infinite dimension. Closed-loop stability is guaranteed under the off-line trained NN-approximation of the feedback gains. 
	Paper \cite{krstic2023neural} extends this framework from first-order hyperbolic PDEs to a more complex class of
	parabolic PDEs whose kernels are governed by second-order PDEs, raising the difficulty for solving such PDEs and for proving the sufficient  smoothness of their solutions, so that the NO approximations have guarantee of sufficient accuracy.
	Furthermore, an operator learning framework for accelerating nonlinear adaptive control is proposed in \cite{bhan2023operator}, where three operators are trained, namely parameter identifier operator, controller gain operator, and control operator.

	\subsection{Results, contributions, and organization of the paper}
	
	In this paper, we employ the DeepONet framework to learn the control kernel functions and the observer gains for the output feedback of a delayed first-order hyperbolic partial integro-differential equation (PIDE) system. Due to the system incorporating state and measurement or actuation delays, two transport PDEs are introduced to represent the delayed states, thus forming a hyperbolic PDEs cascade system. 
	Applying the backstepping transformation, we derive a set of coupled PIDEs that three backstepping kernels should satisfy, the solution of which can only be obtained numerically. Hence, three DeepONets are trained to approximate the three kernel functions from the numerical solutions.  Once the neural operators are trained from data, the kernel equations donot need to be solved numerically again for new functional coefficients and new delays.  
	
	We use the universal approximation theorem to prove the existence of DeepONet approximations, with an arbitrary accuracy, of the exact continuous operators mapping the delay and the system coefficient functions into kernel functions. Based on the approximation result, we provide a  state-feedback stability guarantee under neural operator kernels by using a Lyapunov functional. 
	
	We incorporate a ``dead-time'' into our PIDE model. Dead-time can represent either actuation or sensing delay---when the full state is unmeasured, the delay can be shifted between the input and the output. Without loss of generality, we locate the delay at the output/measurement. In such an architecture, control requires the design of an observer for the unmeasured state. Due to the delayed measurement, the  backstepping transformation for the observer design and analysis contains four kernels, which determine two observer gains. We use two DeepONets to learn the observer gains directly, instead of the four kernel functions. 
	
	The observer with gains produced  by neural  operators is proved to converge to the actual states. Moreover, we prove the stability of the output feedback system under the neural gains through constructing a new Lyapunov functional. Within the proof, we combine the system under prescribed stabilizing controller with the feedback of the estimated states and the observer error system to establish the exponentially stability, thus verifying the separation principle. 
	We demonstrate the theoretical results with numerical tests and the corresponding code is available on \href{https://github.com/JingZhang-JZ/NO_hyperbolic_delay.git}{github}.

	
The paper's main contribution is the following:
	\begin{itemize}
		\item Unlike the two inaugural papers \cite{bhan2023neural, krstic2023neural}, this paper considers a PDE control problem subject to delays. The delayed system brings new challenges of dealing with multi-kernel coupled PDEs, for which we train three NNs to approximate three control kernels and two NNs to approximate two observer gains. 
	\end{itemize}
 
The paper's additional contributions are:
	\begin{itemize}
		\item We combine the observer error system with the closed-loop system under the estimated state feedback to establish the exponentially stability, which verify the separation principle under the DeepONet learned output feedback controller.
		\item We train two DeepONets to approximate the two observer gains  instead of four kernel functions, which cuts the offline training computation cost in half.
	\end{itemize}
	
	The paper is organized as follows.  In Section \ref{Backstepping-method}, we summarize the key steps and related conclusions for designing state feedback, observer design, and output feedback controllers using backstepping methods. Corresponding to the Backstepping theoretical results, Section \ref{DeepONet operator} gives the DeepONet-based design and stability analysis for state feedback, observer and output feedback controller. In particular, we prove the operators are Lipschitz continuous and provide the existence of DeepONet approximations of the operators to any given accuracy in Section \ref{section-Accuracy-approx}. We illustrate the theoretical results with numerical examples in Section \ref{simulation}.
 \color{black} 	
	 
 \textbf{Notation:} Throughout the paper, we adopt the following notations to for functions' domain. 
 \begin{align}\label{define-T1}
 	{\mathcal T}_1&=\left\{(s,q):0\leq  s \leq q \leq 1\right\},\\
 	{\mathcal T}_2&=\left\{(s,r):~0\leq   s,~r  \leq 1 \right\}, \label{define-T2}\\
 	{\mathcal T}_3&=\left\{(s,q):~0\leq   q  \leq s \leq 1, \right\}, \label{define-T3}
 	\\ \label{define-underline C}
 	\underline C&=\left\{c\in C^1[0,1]: ~c(1)=0\right\}.
 \end{align}
 For  $f(s)\in L^p[0,1]$ and $g (s,q)\in L^p (\mathcal{T} )$, where $\mathcal{T}\in \mathbb{R}^2$, with $p=2$ or $\infty$, we define the following norms:
  \begin{align}
     \| f\|:=&\| f\|_\infty =\sup_{s\in [0,1]} |f(s)|, \\
          \|g\|:=&\|g\|_\infty=\sup_{(s,q)\in \mathcal{T}} |g(s,q)|, \\
          \| f\|_{L^2}^2:=&\int_0^1  f^2(s)ds,     \\
    \| g\|_{L^2}^2:=&\int_{\mathcal{T}}  g^2(s,q)dsdq.
  \end{align}

\section{Backstepping Design for a PIDE with Output and State Delays}\label{Backstepping-method}
We consider the following PIDE system with state and sensor delay
  \begin{align}
  	x_t(s,t)=&-x_{s}(s,t) + c(s)x(1,t-\tau)  \nonumber\\&~+ \int^1_{s}   f(s,q)x(q,t)dq,\label{eq:main-x}
  	\\ \label{eq:bnd-x}
  	x(0,t)=&~ U(t),
  	\\\label{eq:measurement-in}
  	y(t)=&~x(1,t-h).
  \end{align}
  for all $(s, t)\in [0,1]\times \mathbb{R}^+$ with the initial condition $x(s, 0)= x_{0}(s)$, and $y(t)$ representing the output that can be measured. There are two types of delays  in the system: recycle delay $\tau$ due to  transport, and measurement delay $h$. The delay $h$ can be alternatively thought of as input delay---the modeler is free to treat ``dead time'' as acting at either the sensor or the actuator. We treat the dead time as acting at the sensor.  
  
\begin{assump}\label{ass-delays}
 Denote the upper bound of the delay by $\bar \tau$, namely, $\tau \le \bar \tau$. Usually, the transportation delay $\tau$ is longer than the dead time $h$, so we  assume $\eta:=\tau-h >0$, and thus $0<h,\eta <\tau\le \bar \tau$. 
\end{assump} 

\begin{assump}\label{ass-1-2}
	$c \in { C}^1 ([0,1])$ with  $c(1)=0$,   $f\in { C}^1( {\cal T}_1)$, and let the following symbols denote their bounds: $\bar c:=\|c\|$, $\bar f :=\|f\|$. 
\end{assump} 

We introduce  transport PDEs to represent the delayed state and delayed measurement, rewriting \eqref{eq:main-x}-\eqref{eq:measurement-in} as:
  \begin{align}
  	x_t(s,t)=&-\!\! x_{s}(s,t)+c(s) u(0,t)+\!\!\!\int_{s}^1 \!\!\!\! f(s,q) x(q,t)dq,  
  	\label{eq:main-x-1}
  	\\ \label{bnd-x}
  	x(0,t)= &~U(t),\\
  	\label{eq:main-v-1}  h v_t(s,t)=&~v_s(s,t), 
  	\\\label{eq:bnd-v-1}v(1,t)=&~ x(1,t),\\
  	\label{eq:main-u-1}
  	\eta u_t(s,t) =&~u_s(s,t),
  	\\\label{eq:bnd-u-1}  u(1,t)=&~v(0,t)=x(1,t-h),
  \end{align}
  for  $(s, t)\in (0, 1) \times \mathbb{R}^+$, with $v_0(s), ~u_0 (s)\in L^2([0, 1])$ denoting the initial conditions for $v$ and $u$, respectively.
We will sketch the backstepping design with state feedback for system \eqref{eq:main-x-1}-\eqref{eq:bnd-u-1} in the following two subsections.

 \subsection{Backstepping design for delay compensation with state-feedback}
  First, we employ the following backstepping transformation: 
  \begin{align}
  	z(s,t)=&~\Gamma|_{K,L,J}[x,v,u](s,t)\nonumber 
  	\\:=&~ x(s,t)-\int^1_s K(s,q)
  	x(q,t)dq\nonumber \\&~ \nonumber
  	-h\int_0^1 L(s+h r )v(r,t)dr\nonumber \\&~-\eta\int_0^1 J(s+\eta r)u(r,t)dr,  \label{eq:trans-ori} 
  \end{align}
  and its associated inverse transformation
  \begin{align}
  	x(s,t)=&~\Gamma^{-1}|_{B,D,E}[z,v,u](s,t)\nonumber\\ :=&~z(s,t)+\int^1_s B(s,q)
  	z(q,t)dq\nonumber \\&~
  	+\int_0^1 D(s,r)v(r,t)dr\nonumber \\&~+\int_0^1 E(s,r)u(r,t)dr,  \label{eq:inv-trans}
  \end{align}
  where kernels $K,~B$ is defined on ${\cal T}_1$, $L$ on $[0,1+h]$ by treating the function $s+hr$ of $(s,r)$  as a single variable, $J$ on $[0,1+\eta]$ by treating the function $s+\eta r$ of $(s,r)$ as a single variable, and $D,~E$ on $ {\cal T}_2$. The task of the transformation \eqref{eq:trans-ori} is to produce the following stable target system:
  \begin{align}\label{eq:main-tar-z}
  	z_t(s,t)=&~-z_{s}(s,t),~\forall (x, t)\in (0, 1) \times \mathbb{R}^+,\\
  	\label{eq:bnd-tar-z}
  	z(0,t)=&~0,\\
  	\label{eq:main-v-t}  h v_t(s,t)=&~v_s(s,t), 
  	\\\label{eq:bnd-v-t}v(1,t)=& ~z(1,t),\\
  	\label{eq:main-u-t}
  	\eta u_t(s,t) =&~u_s(s,t),
  	\\\label{eq:bnd-u-t}  u(1,t)=&~v(0,t).
  \end{align}
  To map \eqref{eq:main-x-1}-\eqref{eq:bnd-u-1} into \eqref{eq:main-tar-z}-\eqref{eq:bnd-u-t}, the kernels need to satisfy:
  \begin{align}\label{eq:K}
  	K_{s}(s,q)=&-K_{q}(s,q) +f(s, q)-\!\!\int^q_s \!\!\!K(s,r)f(r,q)dr,
  	\\\label{eq:K-bnd-2}
  	K(s,1)=  &  L(s+ h),\\
  	\label{eq:L}
  	L(\phi)= &\begin{cases}  
  		J(\phi+\eta),
  		& \phi< 1 \\
  		0, & \phi \ge 1 \\
  	\end{cases},
  	\\\label{eq:J}
  	J(\sigma)= &\begin{cases}  
  		\int_{\sigma}^1 K(\sigma,q)c(q)dq-c(\sigma),
  		& \sigma< 1 \\
  		0, & \sigma \ge 1 \\
  	\end{cases}.
  \end{align}
  As $c(1)=0$ is assumed in Assumption \ref{ass-1-2}, $J$ is continuous at $\sigma=1$.
  Substituting \eqref{eq:J} and \eqref{eq:L} into \eqref{eq:K-bnd-2}, one gets  
  \begin{align}\label{eq:K1-final}
  	&K_{s}(s,q)\!+\!K_{q}(s,q) =
  	f(s, q)\!-\!\int^q_s
  	K(s,r)f(r,q)dr,\\
  	&K(s,1)=\begin{cases}  
  		\int_{s+ \tau}^1   K(s+ \tau,\theta)c(\theta)d\theta
  		-c(s+ \tau), &s+\tau< 1  \\
  		0,&s +\tau \ge 1 \\
  	\end{cases}\label{eq:K-bnd}
  \end{align}
  It is worth  noticing that $K(s,1)=0$ when $\tau \ge 1$, which implies that only one-case situation is needed.
  Using the method of characteristics, we get the integral form
  \begin{align}\label{eq:integral-K}
  	K(s, q)
  	=\begin{cases}(\Phi_0(f)+  \Psi_{0}(c))(s,q)\\  ~+(\Phi(f,K)+ \Psi(c,K))(s,q), &  
  		s+\tau< q \\
  		\Phi_0(f)(s,q)+\Phi(f,K)(s,q), & s+\tau \ge q \\
  	\end{cases}
  \end{align}
  where $\Phi_0$ and $\Psi_0$ are depends on $f$ and $c$, respectively,
  \begin{align}\label{eq:Phi0}
  	\Phi_{0}(f)(s,q)&=-\int^{s+1-q}_{s}
  	f(\theta, \theta-s+q) d\theta,\\
  	\Psi_0(c)(s,q)&=- c(s-q+1+  \tau),\label{eq:Psi0}
  \end{align} 
  and
  $\Phi$ and $\Psi$ are functionals acting on $K$,
  \begin{align}\label{eq:Phi}
  	\Phi(f,K)&= \int^{s+1-q}_{s}\!\!\!\int_{\theta}^{\theta-s+q}
  	f(r,\theta-s+q)K(\theta,r)dr d\theta,
  	\\\label{eq:Psi}
  	\Psi(c,K)&= \int_{s-q+1+ \tau}^1c(\theta)K
  	(s-q+1+ \tau,\theta)d\theta.
  \end{align}
  Based on \eqref{eq:integral-K}, we can derive  $L$, $J$ from  \eqref{eq:L} and \eqref{eq:J}.

\begin{theorem}\label{Th-kernel}
\rm{
For every $(f, c) \in C^1(\mathcal{T}_1) \times \underline C $, the kernel  $K \in C^0 (\mathcal{T}_1)$ and  $L,~J \in C^0 [0,1+h]$ have bounds  	 
  	\begin{align}
  		|K(s,q)|\leq &\left(\bar c+\bar f\right) \mathrm{e}^{ ( \bar c+\bar f)(q-s)}\nonumber \\\leq & \bar K: = \left(\bar c+\bar f\right) \mathrm{e}^{ ( \bar c+\bar f)},\label{bound-K}\\
  		|L(s)| \leq & ~ \bar c \mathrm{e}^{(\bar c+\bar f)(1-s)} \leq \bar L: =\bar c \mathrm{e}^{(\bar c+\bar f)} ,\label{bound-L}\\
  		|J(s)| \leq &\bar L.\label{bound-J}
  	\end{align}}
\end{theorem}

The proof can be found in \cite{qi2021output}. 
  From the boundary conditions \eqref{bnd-x} and \eqref{eq:bnd-tar-z}, the controller is 
  \begin{align}
  	U(t)=&\int_0^1   K(0,q)x(q,t)dq+h\int_0^1   L(hr)v(r,t)dr
  	\nonumber \\&+\eta\int_0^1   J(\eta r)u(r,t)dr. 
  	\label{controller-2}
  \end{align}

\subsection{Backstepping design for the observer and the output-feedback}
 In this subsection, we will briefly  introduce the design of the observer and the output-feedback controller using the backstepping method, and the detailed designed process can be found in \cite{qi2021output}. 
  	The proposed observer is a copy of \eqref{eq:main-x-1}-\eqref{eq:bnd-u-1} with the measurement error:
  	\begin{align}
  		\hat x_t(s,t)=&~-\hat x_{s}(s,t)+\int_{s}^1  f(s,q)\hat x(q,t)dq
  		\label{eq:main-x-o}\\\nonumber
  		& ~+c(s)\hat u (0,t)+Q_{1}(s)(v(0,t)-\hat v(0,t)),~\\ 
  		\hat x(0,t)=&~ U(t),\\
  		\label{eq:main-v-o}   h\hat v_t(s,t)= &~\hat v_s(s,t)
  		+Q_{2}(s)(v(0,t)-\hat v(0,t)),
  		\\\label{eq:bnd-v-o}\hat v(1,t)=&~\hat x(1,t),\\
  		\label{eq:main-u-o}
  		\eta \hat u_t(s,t)
  		= & \hat u_s(s,t),
  		\\\label{eq:bnd-u-o}\hat u(1,t)=& x(1,t-h),
  	\end{align}
  	where  observer gains $Q_1(s),~ Q_2(s) \in  {L}^2(0,1)$ are to be determined later and the initial conditions are denoted by $ \hat x_0, \hat v_{0},\hat u_{0}\in {L}^2(0,1)$.  
  	Define the error states:
  	\begin{align*}
  		\tilde x=x-\hat x, \quad \tilde v=v-\hat v, \quad \tilde u=u-\hat u,
  	\end{align*}
  	which gives
  	\begin{align}
  		\tilde x_t(s,t)=&~-\tilde x_{s}(s,t)+
  		\int_{s}^1  f(s,q)\tilde x(q,t)dq +c(s)\tilde u(0,t)
  		\nonumber \\&~-Q_{1}(s)\tilde v(0,t),\label{eq:main-x-e}\\ \label{eq:bnd-x-e}
  		\tilde x(0,t)=&~0,\\\label{eq:main-u-e}
  		h \tilde  v_t(s,t)=&~  \tilde v_s(s,t)-Q_{2}(s)\tilde v(0,t),
  		\\\label{eq:bnd-u-e}\tilde v(1,t)=&~\tilde x(1,t),
  		\\
  		\eta \tilde u_t (s,t)= &~\tilde u_s (s,t),
  		\\\label{eq:bnd-u2-e} \tilde u(1,t)= &~0, 
  	\end{align}
  	with the initial conditions $\tilde x_0=x_0-\hat x_0$, $\tilde v_0= v_0-\hat v_0$,  $\tilde u_{0}= u_{0}-\hat u_{0}$.
  	We employ the following backstepping transformations,
  	\begin{align}
  		\tilde x(s,t)=&~\mathcal{F}|_{F,M,P}[\tilde z,\tilde w](s,t):=\tilde z(s,t)  \nonumber \\&~ -\int_{s}^{1}  F(s,q)\tilde z(q,t)dq - \int_{0}^{s} M(s,q)\tilde w(q,t)dq\nonumber\\&~- \int^{1}_{s} P(s,q)\tilde w(q,t)dq, 
  		\label{eq:obs-transform-z-x}\\
  		\tilde v(s,t)=&~\mathcal{R}|_{R}[\tilde w](s,t):= \tilde w(s,t) \nonumber\\&~- \int_{0}^{s} R(s-q)\tilde w(q,t) dq,
  		\label{eq:obs-transform-w-v}
  	\end{align}
  	and their associated inverse transformations
  	\begin{align}
  		\tilde z(s,t)= &~\mathcal{F}^{-1}|_{\breve F,\breve M,\breve P}[\tilde x,\tilde v](s,t):=\tilde x(s,t) \nonumber \\&~+ \int_{s}^{1} \breve F(s,q)\tilde x(q,t)dq \nonumber\\&~
  		+\int^{1}_{0}\breve P(s+hq)\tilde v(q,t)dq, 
  		\label{eq:obs-transform-x-z}\\
  		\tilde w (s,t) = &~\mathcal{R}^{-1}|_{\breve R}[\tilde v](s,t):= \tilde v(s,t) \nonumber \\&~+ \int_{0}^{s}\breve R(s-q)\tilde v(q,t)dq, 
  		\label{eq:obs-transform-v-w}
  	\end{align}
  	where observer kernels $M$ defined in $\mathcal{T}_3$,
  	and $F,~\breve{F},~ P$ defined in ${\cal T}_1$, while $\breve{P}$ defined in $[0,1+h]$ by treating the function $s+hq$ of $(s,r)$  as a single variable, and $R,~\breve R$ defined in $[0,1]$ by treating the function $s-q$ of $(s,r)$  as a single variable.
  	The transformations \eqref{eq:obs-transform-z-x} and \eqref{eq:obs-transform-w-v} admit the following observer error target system:
  	\begin{align}
  		\label{eq:target-z}
  		\tilde z_t(s,t)&=-\tilde z_s(s,t)+S(s)\tilde u(0,t),\\
  		\label{eq:target-bnd-z} 
  		\tilde  z(0,t)&=0,\\
  		\label{eq:target-w}
  		h\tilde w_t(s,t)&=\tilde w_s(s,t),
  		\\\label{eq:target-bnd-w}\tilde w(1,t)&=\tilde z(1,t),
  		\\
  		\eta \tilde u_t(s,t)&= \tilde u_s(s,t),
  		\\\label{eq:target-bnd-u2-o}\tilde u (1,t)&= 0,
  	\end{align}
  	where 
   \begin{align}
       S(s)=c(s)+\int_s^1 F(s, q)S(q)dq,\label{Eq:S}
   \end{align}  
   and we have $\bar S:=\|S\|= \bar c \mathrm{e}^{\bar F}$. 
  	To convert the error system to the target system, the observer kernels need to satisfy 
  	\begin{align}
  		\label{eq:main-F1}
  		F_{s}(s,q)=&~-F_{q}(s,q)+
  		\int_s^q f(s,r)F(r,q)dr
  		-f(s,q),\\ \label{eq:main-F2}
  		hM_{s}(s,q)=&~M_{q}(s,q)+  h\int^1_s
  		f(s,r)M(r,q)dr,
  		\\ \label{eq:main-F3}
  		hP_{s}(s,q)
  		=&~P_{ q}(s,q)+h \int^1_q
  		f(s,r)M(r,q)dr
  		\nonumber \\&~  + h\int_s^q
  		f(s,r)P(r,q)dr,
  		\\ \label{eq:bnd-F1}
  		F(0,q)=&~ 0, \quad
  		M(s,s)=P(s,s),\\
  		\label{eq:bnd1-F3}
  		P(0,q)=&~0,\quad
  		P(s,1)=h F(s,1),
  		\\\label{solution-R}
  		R(\xi)=&~M(1,1-\xi),
  	\end{align} 
  	with the observer gains are given 
  	\begin{align}\label{eq:Q1}
  		Q_1(s)=&-\frac{1}{h}M(s,0),\\
  		\label{eq:Q2}
  		Q_2(s)=&-R(s)=-M(1,1-s).
  	\end{align}
  	To realize the inverse transformation, the inverse kernels satisfy 
  	\begin{align}
  		\label{eq:inv-main-F1}
  		\breve F_{s}(s,q)=&-\breve F_{q}(s,q)
  		-\int_s^q f(s,r) \breve F(r,q)dr
  		-f(s,q),
  		\\ \label{eq:inv-bnd-F1}
  		\breve F(0,q)= &~0, \\
  		\label{eq:breve-P}
  		\breve P(\varsigma)= &\begin{cases}  
   h \breve F (\varsigma-h,1)
   & \varsigma > h\\
   0, & \varsigma \le h \\
  		\end{cases},
  		\\\label{inv-solution-R}
  		\breve R(\zeta)=&~\breve P(1+h(1-\zeta)).
  	\end{align}

\begin{theorem}\label{Th-obskernel}
   \rm{
  		For every $(h, f) \in \mathbb{R}^+\times C^1(\mathcal{T}_1) $, the kernel equation \eqref{eq:main-F1}-\eqref{solution-R} admits a unique solution $F, P \in C^0 (\mathcal{T}_1)$, $M\in C^0 (\mathcal{T}_3)$ and $ R\in C^0 [0,1]$ with the bound
  		\begin{align}
   |F(s,q)|\leq &~\bar f  \mathrm{e}^{ \bar f(q-s)} \leq \bar F: = \bar f  \mathrm{e}^{ \bar f} ,\\
   |M(s,q)| \leq &~ h\bar F \mathrm{e}^{2\bar f h(1-q)} \leq \bar M: =h \bar f  \mathrm{e}^{ \bar f (2h+1)} ,\\
   |P(s,q)| \leq &~\bar M, ~~~~~~~|R(\xi)| \leq \bar M.
  		\end{align}
  		Also, the inverse kernel equations \eqref{eq:inv-main-F1}-\eqref{inv-solution-R} admits a unique solution  $\breve F\in C^0 (\mathcal{T}_1)$, $\breve P \in C^0 [0,1+h] $ and $\breve R \in C^0 [0,1]$, with the bound
  		\begin{align}
   |\breve F(s,q)|\leq &  \bar F ,
   \quad |\breve P(\varsigma)| \leq h \bar F, \quad |\breve R(\zeta)| \leq h \bar F.
  		\end{align}
  		Further, the observer gains $Q_1,~Q_2 \in C^0[0,1]$, with bound
  		\begin{align}
   |Q_1(s)|\leq &\bar Q_1:=\bar f  \mathrm{e}^{ \bar f (2h+1)}, ~~~|Q_2(s)|\leq \bar M.
  		\end{align}}
\end{theorem}
The proof is provided in \cite{qi2021output}.
We put together the observer \eqref{eq:main-x-oN}-\eqref{eq:bnd-u-oN}  along with state-feedback controller \eqref{controller-2} and finally obtain the output-feedback controller
  \begin{align}\label{hat-estimated-controll}
  	U(t)=&\int_0^1   K(0,q) \hat x(q,t)dq+h\int_0^1    L(hr)\hat v(r,t)dr
  	\nonumber \\&+\eta \int_0^1  J(\eta r)\hat u(r,t)dr,
  \end{align}
 which stabilizes the system \eqref{eq:main-x-1}-\eqref{eq:bnd-u-1}.

\section{Output and State Delay Compensation under DeepONet based Controller}\label{DeepONet operator}
Before proceeding, we first present the following theorem on the DeepONet approximability of operators between function spaces.

  \begin{theorem}\label{Th-deeponet}\rm{
  	(DeepONet universal approximation theorem \cite{deng2022approximation}, Theorem 2.1). Let $X\in \mathbb{R}^{d_x}$ and $Y\in \mathbb{R}^{d_y}$ be compact sets of vectors $x\in X$ and $y\in Y$, respectively. Let $\mathcal{U}: X\rightarrow U \subset \mathbb{R}^{d_u}$ and $\mathcal{V}: Y\rightarrow V \subset \mathbb{R}^{d_v}$ be sets of continuous functions $u(x)$ and $v(y)$, respectively. Let $\mathcal{U}$ be also compact. Assume the operator $\mathcal{G}: \mathcal{U}\rightarrow\cal{V} $ is continuous. Then for all $\varepsilon>0$, there exist $m^*$, $p^*\in \mathbb{N}$ such that for each $m\ge m^*$, $p\ge p^*$, there exist $\theta^{(k )}$, $\vartheta ^{(k )}$ for neural networks $f^{\cal{N}}(\cdot;\theta^{(k )})$,  $g^{\cal{N}}(\cdot;\vartheta^{(k )})$, $k=1,...,p,$ and $x_j \in X$, $j=1,...,m,$ with corresponding $\mathbf{u}_m=(u(x_1), u(x_2),...,u(x_m))^T$, such that
  	\begin{equation}
  		|\mathcal{G}(u)(y)-\mathcal{G}_\mathbb{N}(\mathbf{u}_m)(y)|\leq \varepsilon,
  	\end{equation}
  	where 
  	\begin{equation}
  		\mathcal{G}_\mathbb{N}(\mathbf{u}_m)(y)=\sum_{k=1}^{p}g^{\cal{N}}(\mathbf{u};\vartheta^{(k)})f^{\cal{N}}(y;\theta^{(k)}), 
  	\end{equation}
  	for all functions $u\in \mathcal{U}$ and for all values $y\in Y$ of $\mathcal{G}(u) \in \mathcal{V}$. }
  	\end {theorem}
  	
  	The theorem provides the theoretical underpinning for the utilization of DeepONet-based controllers, enabling the approximation of control kernel operators using neural networks if the they are continuous.   
  	In this section, we will utilize three DeepONet to approximate the three state-feedback control kernel operators and two DeepONets to approximate two observer gain operators, instead of four observer kernel operators. These operators are defined as follows:

\begin{definition}
  		Kernel operator $\mathcal{K}: \mathbb{R}^+\times C^1(\mathcal{T}_1)\times \underline C \mapsto C^0(\mathcal{T}_1)  $, $\mathcal{L}:\mathbb{R}^+\times \mathbb{R}^+\times  C^1(\mathcal{T}_1)\times \underline C \mapsto C^0[0,1+h]  $ and 
  		$\mathcal{J}:\mathbb{R}^+\times C^1(\mathcal{T}_1)\times \underline C \mapsto C^0[0,1+\eta]  $
  		are defined by  
  		\begin{align}
  			K(s,q):&=\mathcal{K}( \tau, f, c),\label{ope-K}\\
  			L(\phi):&=\mathcal{L}( \tau, \eta, f, c),\label{ope-L}\\
  			J(\sigma):&=\mathcal{J}( \tau, f, c).\label{ope-J}
  		\end{align}   		
  		For each constant $\tau,~h\in \mathbb{R}^+$ and function $f \in C^1(\mathcal{T}_1)$, $c\in \underline C$, the operators $\mathcal{K}$, $\mathcal{L}$ and $\mathcal{J}$ can generate the kernel functions $K(s,q)$, $L(\phi)$ and $J(\sigma)$, which satisfy equations \eqref{eq:K}-\eqref{eq:J}. 
\end{definition}

It is worth noting that the operator $\mathcal{L}$ is independent of $h$ because $h$ solely affects the domain of $\phi$, as shown in equation \eqref{eq:L} that $\mathcal{L}(\phi)=0$ if $\phi\ge 1$. Similarly, the operator $\mathcal{J}$ is independent of $h$.

\begin{definition}
  	Observer gain $\mathcal{Q}_1,~ \mathcal{Q}_2: \mathbb{R}^+\times  C^1(\mathcal{T}_1) \mapsto C^0[0,1]  $ are defined by 
  	\begin{align}\label{ope-Qi}
  		Q_i(s):&=\mathcal{Q}_i( h, f),~~i=1,2
  	\end{align}
  	where ${\cal{T}}_1$ is defined in \eqref{define-T1}.
\end{definition}
  	It is noteworthy that we employ directly NNs to train the operators for the observer gains. This choice is driven by both the considerable number of observer kernels (four in total) and the fact that only the gains play a role in the observer's functioning.

\subsection{Accuracy of approximation of backstepping operator with DeepONet}\label{section-Accuracy-approx}
\begin{lemma}\label{lem-Lipschitz-kernel}
(Lipschitzness of backstepping kernel operators). The kernel operators $\mathcal{K}: \mathbb{R}^+\times C^1(\mathcal{T}_1)\times \underline C \mapsto C^0(\mathcal{T}_1)  $, $\mathcal{L}:\mathbb{R}^+\times \mathbb{R}^+\times  C^1(\mathcal{T}_1)\times \underline C \mapsto C^0[0,1+h]  $ and 
  		$\mathcal{J}: \mathbb{R}^+\times C^1(\mathcal{T}_1)\times \underline C \mapsto C^0[0,1+\eta]  $ are locally Lipschitz and, specifically, for any $\bar \tau, ~\bar f, \bar c$, the operators satisfy
  		\begin{align}
  			& \|\mathcal{K}(\tau_1, f_1, c_1)-\mathcal{K}(\tau_2, f_2, c_2)\| \label{K-Lip} \\
  			 \leq &L_k \max\{|\tau_1-\tau_2|, \|f_1-f_2\|, \|c_1-c_2\|\},\nonumber\\
  			 & \|\mathcal{L}(\tau_1, \eta_1, f_1, c_1)-\mathcal{L}(\tau_2,\eta_2, f_2, c_2)\|  \label{L-Lip}  \\
  			 \leq &L_L \max\{|\tau_1-\tau_2|, |\eta_1-\eta_2|, \|f_1-f_2\|, \|c_1-c_2\|\},\nonumber\\
  			 & \|\mathcal{J}(\tau_1, f_1, c_1)-\mathcal{J}(\tau_2, f_2, c_2)\|  \label{J-Lip} \\
  			 \leq &L_J \max\{|\tau_1-\tau_2|, \|f_1-f_2\|, \|c_1-c_2\|\}, \nonumber
  		\end{align}
  	with the Lipschitz constant  $L_K, L_L, L_J>0$.
\end{lemma}

\begin{proof}
We begin with the Lipschitz continuity of operator $\mathcal{K}$, rewriting left hand side of \eqref{K-Lip} as 
\begin{align}
    \|K\tau_1-K\tau_2+Kf_1-Kf_2+Kc_1-Kc_2\| ,
\end{align}
where
\begin{eqnarray}
    &K\tau_1-K\tau_2= \mathcal{K}(\tau_1, f_1, c_1)-\mathcal{K}(\tau_2, f_1, c_1),\\
   & Kf_1-Kf_2= \mathcal{K}(\tau_2, f_1, c_1)-\mathcal{K}(\tau_2, f_2, c_1), 
    \\
    & Kc_1-Kc_2= \mathcal{K}(\tau_2, f_2, c_1)-\mathcal{K}(\tau_2, f_2, c_2). \label{Kc1-Kc2}
\end{eqnarray}
We first consider the continuity of $\mathcal{K}$ w.r.t. $\tau$, as $\tau$ is a scalar parameter. 
The integration form \eqref{eq:integral-K}, can be further rewritten in the term of operator with two branches: 
\begin{align}
    \mathcal{K}_1(\tau) =&\Phi_0(f) +\Psi_0(c,\tau)+\Phi_{11}(\tau, f, \mathcal{K}_1(\tau)) \nonumber \\
    &+\Phi_{12}(\tau, f, \mathcal{K}_2) +\Psi_{11}(\tau, c, \mathcal{K}_1(\tau)) \nonumber \\
    & +\Phi_{12}(\tau, c, \mathcal{K}_2),~~ \textit{for}~s+\tau>q, \label{eq:K1-tau0}
    \\
    \mathcal{K}_2  =&\Phi_0(f)  +\Phi(f, \mathcal{K}_2(\tau)),~ ~\textit{for}~s+\tau\le q .
    \label{eq:K2-notau}
\end{align}
where
\begin{align}
    &\Phi_{11}(\tau,f,\mathcal{K}_1(\tau))(s,q)\nonumber \\
    &~= \int^{s+1-q}_{s} 
    \!\!\!\int_{\theta+\tau}^{\theta-s+q}
  	f(r,\theta-s+q)\mathcal{K}_1(\tau)(\theta,r)dr d\theta,\label{Phi11}\\
   &\Phi_{12}(\tau,f,\mathcal{K}_2)(s,q)\nonumber \\
    &~= \int^{s+1-q}_{s} 
    \!\!\!\int_{\theta}^{\theta+\tau}
  	f(r,\theta-s+q)\mathcal{K}_2(\theta,r)dr d\theta,\label{Phi12}\\
 &\Psi_{11}(\tau,f,\mathcal{K}_1(\tau))(s,q)\nonumber \\
    &~= \int^{1}_{\psi(\tau,s,q)}   
  	c(\theta)\mathcal{K}_1(\tau)(s-q+1+\tau,\theta)d\theta,\label{Psi11}\\
   &\Psi_{12}(\tau,f,\mathcal{K}_2)(s,q)\nonumber \\
    &~= \int_{s-q+1+\tau}^{\psi(\tau,s,q)}      
  	c(\theta)\mathcal{K}_2(s-q+1+\tau,\theta)d\theta \label{Psi12},
  \end{align}
  with 
  \begin{align}
      \psi(\tau,s,q)=\min\{1,s-q+1+2\tau\}.
  \end{align}
Take the derivative of the operators of $\tau$, 
\begin{align}
    \partial_\tau  \mathcal{K}_1(\tau) =&\Gamma(s,q)+ \Phi_{11}(\tau, f, \partial_\tau\mathcal{K}_1(\tau))\nonumber 
    \\&+\Psi_{11}(\tau, c, \partial_\tau \mathcal{K}_1(\tau)),~~ \textit{for}~~s+\tau>q, \label{eq:K1-tau}
    \\
   \partial_\tau  \mathcal{K}_2  =& 0,~ ~~~~~~ ~~ ~~~~~~ ~~~~~~~~~~~ ~~~~~~ ~~\textit{for}~~s+\tau\le q.
\end{align}
where 
\begin{align}
    \Gamma&(s,q)=c'(s-q+1+\tau) \label{Gamma-1} \\
    &+\int_s^{s+1-q}f(\theta+\tau,\theta-s+q)[\mathcal{K}_2-\mathcal{K}_1(\tau)](\theta,\theta+\tau)d\theta \nonumber\\
    &+2c(\psi(\tau,s,q))[\mathcal{K}_2-  \mathcal{K}_1(\tau)](s-q+1+\tau,\psi(s,q,\tau))\nonumber \\
    &-c(s-q+1+\tau)\mathcal{K}_2(s-q+1+\tau,s-q+1+\tau)
    \nonumber \\&+
     \Psi_{11}\left(\tau, c,  \partial_s\left(\mathcal{K}_1(\tau)(s,\theta)\right) \right)  +\Phi_{12}\left(\tau, c, \partial_s\left(\mathcal{K}_2 (s,\theta)\right)\right).\nonumber
\end{align}
Notice that $\partial_s\left(\mathcal{K}(s,\theta)\right)=K_s(s,\theta)$. It is straightforward to demonstrate in a similar way the proof of Theorem \ref{Th-kernel} that $K_s(s,q)$ is bounded by 
\begin{equation}
     |K_s(s,q)|\le \Gamma_0 \mathrm{e}^{\bar c (q-s)}, \label{eq:bound_Ks}
\end{equation} 
  with a constant $\Gamma_0>0$. Addition to the fact $K$ is bounded, $\Gamma $ is also bounded and the bounds is denoted by $\bar \Gamma =:\|\Gamma(s,q)\| $.
Applying the successive approximation approach, we reach the boundeness of $\partial_\tau\mathcal{K}_1(\tau)$ as follows
\begin{align}
     |\partial_\tau\mathcal{K}(\tau)|\leq &\bar \Gamma \sum_{n=0}^\infty
     (\bar c+\bar f)^{n+1}\frac{(q-s)^n}{n!}\nonumber\\
     = &\bar \Gamma(\bar c+\bar f)  \mathrm{e}^{(\bar c+\bar f)(q-s)}
     \leq \bar \Gamma(\bar c+\bar f)  \mathrm{e}^{(\bar c+\bar f)}.
\end{align}
Consequently, we infer that operator $\mathcal{K}$ is Lipschitz continuous of $\tau$ with Lipschitz constant $\bar \Gamma(\bar c+\bar f)  \mathrm{e}^{(\bar c+\bar f)}$. 
Second, we investigate the boundedness of $Kf_1-Kf_2$.
From \eqref{eq:integral-K}, we have
\begin{align}
    Kf_1-Kf_2= & \Phi_0(f_1-f_2)+\Phi(f_1-f_2,\mathcal{K}(f_2))\nonumber \\&+\Phi(f_1,Kf_1-Kf_2),
\end{align}
Introduce the iteration
\begin{align}
    \delta_f K^{n+1}&=\Phi(f_1,\delta_f K^n),\\
    \delta_fK^0&=\Phi_0(f_1-f_2)+\Phi(f_1-f_2,\mathcal{K}(f_2)),
\end{align}
which verifies 
\begin{align}
    Kf_1-Kf_2=&\sum_{n=0}^\infty \delta_f K^n.
\end{align}
Recalling $K$ is bounded and combining the definition of $\Phi_0$ and $\Phi$ in \eqref{eq:Phi0} and \eqref{eq:Phi}, respectively, we get
\begin{align}
    \|\delta_fK^0\| &= (1+\|K\|)\|f_1-f_2\|.
\end{align}
By induction, 
\begin{align}
    \delta_f K^{n}\leq & (1+\|K\|) \frac{\bar f^n (q-s)^n}{n!}\|f_1-f_2\|.
\end{align}
Therefore it follows that for all $(s,q)\in \mathcal{T}_1$,
\begin{align}
    |Kf_1-Kf_2| \le& (1+\|K\|)  \mathrm{e}^{\bar f}\|f_1-f_2\|.
\end{align}
Third, we consider the boundedness of \eqref{Kc1-Kc2}. From \eqref{eq:integral-K}, it derives
\begin{align}
&Kc_1-Kc_2\\
\nonumber &
  	~~=\begin{cases}c_1-c_2+\Psi\left(c_1\mathcal{K}(c_1)-c_2\mathcal{K}(c_2)\right),&  
  		s+\tau< q\\
  		0, & s+\tau \ge q \\
  	\end{cases}
\end{align}

In a similar way to get the bound of $Kf_1-Kf_2$, it arrives at
\begin{align}
    |Kc_1-Kc_2|\leq & (1+\|K\|) \mathrm{e}^{\bar c(q-s)}\|c_1-c_2\|\nonumber 
    \\ \le& (1+\|K\|)  \mathrm{e}^{\bar c}\|c_1-c_2\|.
\end{align}
Consequently, we get the Lipschitz constant 
\begin{align}
    L_K=3\max\{\bar \Gamma(\bar c+\bar f),(1+(\bar c+\bar f)\mathrm{e}^{(\bar c+\bar f)})\}\mathrm{e}^{(\bar c+\bar f)}.
\end{align}
So far we have shown that  operator $\mathcal{K}$ exhibits local Lipschitz continuity with respect to  inputs $\tau$, $f$ and $c$.

Next, we will prove that operators $\mathcal{L}$ and $\mathcal{J}$ are local Lipshcitz continuous using the similar approach due to they dependent on $\mathcal{K}$ as shown in \eqref{eq:L} and \eqref{eq:J}. Subsequently, we just present the differences from the above proof. Denote 
 \begin{align*}
    & \Theta(c,K)(\sigma)=\int_\sigma^1 \mathcal{K}(\tau,f,c)(\sigma,q)c(q)dq,\\
    & K_1=\mathcal{K}(\tau_1,f_1,c_1),~~~~~K_2=\mathcal{K}(\tau_2,f_2,c_2),
 \end{align*}
and thus 
\begin{align}
    & \|\mathcal{J}(\tau_1, f_1, c_1)-\mathcal{J}(\tau_2, f_2, c_2) \|\nonumber 
    \\   
=&\| \Theta\left(c_1-c_2,K_1\right) +\Theta\left(c_2,K_1-K_2\right)+c_2(\sigma)-c_1(\sigma)\|\nonumber 
    \\\leq &\bar c L_K/3  (|\tau_1-\tau_2|+\|f_1-f_2\|
    +\|c_1-c_2\|)
    \nonumber \\
    & +(1+\|K\|)\|c_1-c_2\|
    \nonumber \\   
    \leq  & L_J/3  (|\tau_1-\tau_2|+\|f_1-f_2\|+\|c_1-c_2\|),
\end{align}
where
\begin{align}
    L_J=\bar c L_K+3\|K\|+3.
\end{align}
Since $L$ is a shift of $J$, we have
\begin{align}
    & \mathcal{L}(\tau_1, \eta_1,f_1, c_1)-\mathcal{L}(\tau_2,\eta_1, f_2, c_2) 
    \nonumber \\
     =& \mathcal{J}(\tau_1, f_1, c_1)(\phi+\eta_1)-\mathcal{J}(\tau_2, f_2, c_2)(\phi+\eta_2) 
    \nonumber \\
     =& \int_{\phi+\eta_1}^1[K_1(\phi+\eta_1,q)-K_2(\phi+\eta_2,q)]c_1(q)dq
     \nonumber \\
     &+\int_{\phi+\eta_1}^{\phi+\eta_2} K_1(\phi+\eta_2,q) c_1(q)dq-c_1(\phi+\eta_1)+c_1(\phi+\eta_2)
    \nonumber \\
    &+\Theta(c_1-c_2,K_1)(\phi+\eta_2)+\Theta(c_2,K_1-K_2)(\phi+\eta_2)
    \nonumber \\
    &-c_1(\phi+\eta_2)+c_2(\phi+\eta_2).
   \end{align}
Recalling $\|K(s_1,q)-K(s_2,q)\|\le \Gamma_0 \mathrm{e}^{\bar c }$ due to \eqref{eq:bound_Ks}, the left hand side of \eqref{L-Lip} becomes
\begin{align}
    & \|\mathcal{L}(\tau_1, \eta_1,f_1, c_1)-\mathcal{L}(\tau_2,\eta_1, f_2, c_2) \|
    \nonumber \\
    \leq &\bar c \Gamma_0 \mathrm{e}^{\bar c} |\eta_1-\eta_2|+\bar c\|K\||\eta_1-\eta_2|+L_c|\eta_1-\eta_2|
    \|
    \nonumber \\
      &+\bar c L_k/3 (|\tau_1-\tau_2|+\|f_1-f_2\|
    +\|c_1-c_2\|) 
    \nonumber \\
    &+(1+\|K\|)\|c_1-c_2\|
     \nonumber \\
    \le & L_L/4(|\tau_1-\tau_2|+\eta_1-\eta_2|+\|f_1-f_2\|+\|c_1-c_2\|), \nonumber 
\end{align}
where 
\begin{align}
    L_L=4\max\{& \bar cL_K/3, 1+\|K\|+\bar c L_K/3, \nonumber \\&\bar c \Gamma_0 \mathrm{e}^{\bar c}+\bar c\|K\|+L_c\},
    \end{align}
 with   $L_c$ is the Lipschitz constant for function $c\in C^1$.
 \end{proof}
  	
   \begin{lemma}\label{lem-Lipschitz-gain}
  	(Lipschitzness of observer gain operators). The observer gain operators $\mathcal{Q}_i: \mathbb{R}^+\times C^1(\mathcal{T}_1) \mapsto C^0[0,1]$ for $i=1,2$, 
  	are locally Lipschitz and, specifically, the operators satisfy
  	\begin{align}\label{Ineq-DQi}
  		 &\|\mathcal{Q}_i(h_1 f_1)-\mathcal{Q}_i(h_2, f_2)\|   
  		\\ \leq&  L_{Qi} \max\{|h_1-h_2|, \|f_1-f_2\|\}, 
  	\end{align}
  	with the Lipschitz constants $L_{Qi}>0$. 
\end{lemma}
The proof of this lemma is similar to that of Lemma \ref{lem-Lipschitz-kernel}, so we omit it due to space constraints.

\color{black}
Based on Theorem \ref{Th-kernel}, Theorem\ref{Th-deeponet}, Lemma \ref{lem-Lipschitz-kernel} and the Theorem 3.3 in paper \cite{deng2021convergence}, we get the following result for the approximation of the kernels by DeepONets. 

\begin{theorem} \label{Th-Nkernel}
   \rm{
  		For any  $( \tau, f, c) \in \mathbb{R}^+\times C^1(\mathcal{T}_1) \times \underline C $ with $\tau<\bar \tau$, $\|f\|=\bar f$ and $\|c\|=\bar c$ and $\varepsilon >0$, there exist positive integers $p^*(\varepsilon)$, $m^*(\varepsilon)$, such that for any $p>p^*$ and $m>m^*$, there are neural networks  $f_i^{\cal{N}}(\cdot;\theta_{i}^{(k )})$,  $g_i^{\cal{N}}(\cdot;\vartheta_{i}^{(k)})$, $i=1,2,3$, $k=1,...,p,$ and $(s,q)_j \in \mathcal{T}_1$, $j=1,...,m$, such that
\begin{align}
  	 &~   |\mathcal{K}(\tau,f,c)-\hat{ \mathcal{K} }((\tau,f,c)_{m})(s,q)|\leq \varepsilon,\nonumber \\
         \label{DON-K} &\hat {\mathcal{K}}   =\sum_{k=1}^{p}g^{\cal{N}}_1((\tau,f,c)_{m};\vartheta_1^{(k)})f^{\cal{N}}_1((s,q);\theta_1^{(k)}), \\
      & ~|\mathcal{L}(\tau,\eta,f,c)-\hat{ \mathcal{L} }((\tau,\eta,f,c)_{m})|\leq \varepsilon,\nonumber \\
       \label{DON-L} & \hat {\mathcal{L}}  =\sum_{k=1}^{p}g^{\cal{N}}_2((\tau,f,c)_{m};\vartheta_2^{(k)})f^{\cal{N}}_2((\phi);\theta_2^{(k)}), \\
    &~   |\mathcal{J} (\tau,f,c)-\hat {\mathcal{J}} ((\tau,f,c)_{m})|\leq \varepsilon,\nonumber \\        
         \label{DON-J} &\hat {\mathcal{J}}  =\sum_{k=1}^{p}g^{\cal{N}}_3((\tau,f,c)_{m};\vartheta_3^{(k)})f^{\cal{N}}_3((\sigma);\theta_3^{(k)}),  
  \end{align}  	
   holds, where $(\tau,f,c)_{m}$ are $3\times m$ matrix, with each row containing $m$ discretized elements for constant $\tau$ and functions $f(s,q)$, $c(s)$. Similarly, ${\tau,\eta, f,c}_{m}$ are $4\times m$ matrix.
   
Furthermore, given a required error tolerance $\varepsilon > 0$, the DeepONets defined in \eqref{DON-K}-\eqref{DON-J} approximate the local Lipschitz continuous kernel operators defined in \eqref{ope-K}-\eqref{ope-J} by employing the number of data point evaluations for $\tau$, $\eta$, $f(s,q)$ and $c(s)$, respectively, on the order of 
\begin{equation}
     m^* \sim \varepsilon^{-1}, 
\end{equation}
the number of basis components in the interpolation when reconstructing the kernel functions spaces on the order of
\begin{align}
   p^* \sim \varepsilon^{-\frac{1}{2}}.
\end{align}}
\end{theorem}

\begin{rmk}
It is worth noting that the parameter $m$ determines the number of grids used for discretizing the function. For instance, a two-dimensional function $f$ should be discretized on a grid on  $\mathcal{T}_1$  with $m$ grid points.
\end{rmk}

\begin{theorem} \label{Th-NObs}
   \rm{
  		For any  $(h,f) \in \mathbb{R}^+\times C^1(\mathcal{T}_1) $  and $\varepsilon >0$, there exist positive integers $p^*(\varepsilon)$, $m^*(\varepsilon)$, such that for any $p>p^*$ and $m>m^*$, there are neural networks  $f_{oi}^{\cal{N}}(\cdot;\theta^{(k )})$,  $g_{oi}^{\cal{N}}(\cdot;\vartheta^{(k)})$, $i=1,2$, $k=1,...,p,$ and $(s,q)_j \in \mathcal{T}_1$, $j=1,...,m$, such that
\begin{align}
&  |\mathcal{Q}_i(h,f)-\hat{ \mathcal{Q} }_i((h,f)_{m})(s)|\leq \varepsilon,\nonumber \\
         \label{DON-Qi} &\hat {\mathcal{Q}}_i=\sum_{k=1}^{p}g^{\cal{N}}_{oi}
    ((h,f)_{m};\vartheta_{oi}^{(k)})f^{\cal{N}}_{oi}(s;\theta_{oi}^{(k)}), 
        \end{align}  	
   holds, where $(h,f)_{m}$ are $2\times m$ matrix, with each row containing $m$ discretized elements for constant $h$ and function $f(s,q)$. 

Furthermore, given a required error tolerance $\varepsilon > 0$, the DeepONets defined in \eqref{DON-Qi} approximate the local Lipschitz continuous operator of the observer gains defined in \eqref{ope-Qi} by employing the number of data point evaluations for $h$ and $f(s,q)$, respectively, on the order of 
\begin{align}
    m^* \sim \varepsilon^{-1}, 
\end{align} 
the number of basis components in the interpolation when reconstructing the kernel functions spaces on the order of
\begin{align}
   p^* \sim \varepsilon^{-\frac{1}{2}}.
\end{align}}
\end{theorem}
\color{black}

 \subsection{State-Feedback Stabilization under DeepONet Gain}
\label{Stability-feedback}
Let $\hat K =: \hat{\mathcal{K}} ( \tau, f,c)$, $\hat L =: \hat{\mathcal{L}} ( \tau,h, f,c)$ and  $\hat J =: \hat{\mathcal{L}} ( \tau,h, f,c)$ be approximate operators, and their image functions, with accuracy $\varepsilon$ relative to the exact backstepping kernel $ K= {\mathcal{K}}( \tau, f,c)$,  $ L= {\mathcal{L}}(\tau, h, f,c)$  and $J= {\mathcal{L}}(\tau, h, f,c)$, respectively. 
The  following theorem establishes the properties of the feedback system. 
   
\begin{theorem}\label{th-fb-st}
   \rm{
  		For any  $(\tau, h, f, c) \in \mathbb{R}^+ \times \mathbb{R}^+ \times C^1(\mathcal{T}_1) \times \underline C $, there exist a sufficiently small $\varepsilon^*>0$, such that the feedback control law 
  		\begin{align}\label{hat-controller}
   U(t)=&\int_0^1  \hat K(0,q)x(q,t)dq+h\int_0^1   \hat L(hr)v(r,t)dr
   \nonumber \\&+\eta \int_0^1 \hat J(\eta r)u(r,t)dr .
  		\end{align}
  		with NO gain kernel $\hat K$ and its derived kernels $\hat L$ and $\hat J$ of approximation accuracy $\varepsilon \in (0,\varepsilon^*)$ ensures that the closed-loop system satisfies the exponential stability bound, for all $t>0$
  		\begin{align}\label{Th-feedback-ineq}
   &\|x\|^2_{L^2}+\|v\|^2_{L^2} +\|u\|^2_{L^2}\\
   \leq &~W_0 \mathrm{e}^{-\alpha_0t} (\|x_0\|^2_{L^2}+\|v_0\|^2_{L^2}+\|u_0\|^2_{L^2}),  
  		\end{align}
  		with $W_0>0$ and $\alpha_0>0$. 
    }
\end{theorem}

\begin{proof}
Before proceeding, let $\tilde K =K -\hat K $, $\tilde L=L-\hat L$ and $\tilde J=J-\hat J$ denote the difference between the kernels and their  approximations.  
  		
The proof includes three steps. First, we take the same transformation as \eqref{eq:trans-ori},
  		while with the controller \eqref{hat-controller}, we have the following target system:
  		\begin{align}
   z_t(s,t)=& ~- z_{s}(s,t),\label{hat-target-z}\\ 
    z(0,t)=& -\int_0^1 \tilde K(0,q)x(q,t)dq-h\int_0^1 \tilde L(hr)v(r,t)dr\nonumber
   \\&-\eta\int_0^1\tilde J(\eta r)u(r,t)dr,  \label{hat-bnd-z}
   \\ \label{eq:main-tar-v}
   h v_t(r,t)=& ~v_s(r,t), 
   \\\label{eq:bnd-tar-v}
   v(1,t)=& ~z(1,t),\\
   \label{eq:main-u1}
   \eta u_t(r,t)=& ~u_s(r,t),\\\label{eq:bnd-u1}
   u(1,t)=& ~v(0,t).
\end{align}
Second, we substitute the inverse transformation of \eqref{eq:trans-ori} into \eqref{hat-bnd-z} and get a boundary condition exclusively containing states $(z,v,u)$
\begin{align}\label{N-z0}
    z(0,t)=&-\int_0^1 \tilde K(0,q)\Gamma^{-1}|_{B,D,E}[z,v,u](q,t)dq
    \\&-h\int_0^1 \tilde L(hr)v(r,t)dr-\eta\int_0^1\tilde J(\eta r)u(r,t)dr\nonumber 
    \\=&-\int_0^1 z(q,t)[\tilde K(0,q)+\int_0^q \tilde K(0,r)B(r,q)dr]dq\nonumber 
    \\&-\int_0^1 v(q,t)[h\tilde L(hq)+\int_0^1 \tilde K(0,r)D(r,q)dr]dq\nonumber 
    \\&-\int_0^1 u(q,t)[\eta\tilde J(\eta q)+\int_0^1 \tilde K(0,r)E(r,q)dr]dq.\nonumber 
\end{align}
Substituting \eqref{eq:inv-trans} into \eqref{eq:trans-ori}, we get the relationship between the direct and inverse backstepping kernels:
  		\begin{align}
   B(s,q)&= K(s,q)+\int_s^q  K(s,a)(a,q)da,\\
   D(s,r)&=h L(s+hr)+\int_s^1  K(s,a) D(a,s)da,\\
    E(s,r)&=\eta  J(s+\eta r)+\int_s^1  K(s,a) E(a,s)da.
  		\end{align}
Hence, the inverse kernel satisfies the following bounds:
\begin{align}
   \|B\|\leq&\bar B:= \|  K\| \mathrm{e}^{\| K\|}=\bar K\mathrm{e}^{\bar K},\label{hat-B-bound}\\
   \| D\|\leq& \bar D:=   h\|  L\| \mathrm{e}^{\|  K\|}=h\bar L\mathrm{e}^{\bar K},\label{hat-D-bound} \\
   \|  E\|\leq& \bar E:= \eta \| J\| \mathrm{e}^{\| K\|}=\eta\bar L\mathrm{e}^{\bar K}.\label{hat-E-bound}
\end{align}
Third, we carry out the Lyapunov stability analysis. Define the following Lyapunov functionals:
  		\begin{align} \label{def-V12}
   V_1&=\|z\|_{L^2}^2,~~~
   V_2  =\int_0^1 \mathrm{e}^{-b_1 s}z^2(s,t)ds, \\
   \label{def-V34}
   V_3&=\| v\|_{L^2}^2,~~~
   V_4 =h\int_0^1 \mathrm{e}^{b_2s} v^2(s,t)ds,\\
   \label{def-V56}
   V_5&=\| u\|_{L^2}^2,~~~
   V_6 =\eta \int_0^1 \mathrm{e}^{b_3s} u^2(s,t)ds,
  		\end{align}
  		with $b_i>0$, $i=1,2,3$.
  		Note that the following Lyapunov functional pairs satisfy norm-equivalence relationships: $V_1$ and $V_2$;  $V_3$ and $V_4$; $V_5$ and $V_6$, namely, 
  \begin{align}\label{norm-equt1}
 V_2&\leq V_1\leq \mathrm{e}^{b_1}V_2,\\
 \frac{1}{h}\mathrm{e}^{-b_2} V_4&\leq V_3 \leq \frac{1}{h}V_4,\label{norm-equt2}\\
  \frac{1}{\eta}\mathrm{e}^{-b_3} V_6&\leq V_5 \leq \frac{1}{\eta}V_6.\label{norm-equt3}
  \end{align}
Taking the time derivative of  $V:=\beta_1 V_2+\beta_2 V_4 +V_6$ with $\beta_i>0$, $i=1,2$, we have 
\begin{align*}
   \dot V=&- \beta_1 \int_0^1 -\mathrm{e}^{-b_1s}(z^2(s))_s ds+ \beta_2\int_0^1 -\mathrm{e}^{b_2s}(v^2(s))_s ds 
   \\&+ \int_0^1 -\mathrm{e}^{b_3s}(u^2(s))_s ds  \nonumber\\
   =& -\beta_1 \mathrm{e}^{-b_1} z^2(1)+\beta_1 z^2(0)-\beta_1 b_1 V_2+\beta_2  \mathrm{e}^{b_2} v^2(1)\nonumber \\
   &-\beta_2 v^2(0)-\beta_2 b_2  V_4/h+ \mathrm{e}^{b_3} u^2(1)-u^2(0)-b_3 V_6/\eta.\nonumber    
\end{align*}
Recalling the boundary condition \eqref{N-z0} and $|\tilde K|, |\tilde L|, |\tilde J| <\varepsilon$ given in  Theorem \ref{Th-Nkernel}, we know
\begin{align}
z^2(0,t)\leq &6\varepsilon^2 \breve K V_1+6\varepsilon^2 h^2  \breve L V_3
+
6\varepsilon^2 \eta^2 \breve L V_5,\label{ieq:z0}
\end{align}
where 
\begin{align}\label{define-breve_K_L}
    \breve K=1+\bar K^2 \mathrm{e}^{2\bar K},~~~
\breve L=1+\bar L^2 \mathrm{e}^{2\bar K}.
\end{align}
In addition to the norm inequalities \eqref{norm-equt1}-\eqref{norm-equt3}, we reach
\begin{align}\label{derivative-V1}
   \dot V\leq 
   & -(\beta_1 \mathrm{e}^{-b_1}-\beta_2 \mathrm{e}^{b_2}) z^2(1)
   -(\beta_2- \mathrm{e}^{b_3})v^2(0)-u^2(0)\nonumber 
   \\&-(b_1\beta_1 -6\varepsilon^2 \beta_1\mathrm{e}^{b_1} \breve K)V_2 
   -(b_2\beta_2/h -6\varepsilon^2\beta_1 h \breve L)V_4\nonumber 
   \\&
   -(b_3/\eta  -6\varepsilon^2\beta_1\eta \breve L)V_6.    
\end{align}
Letting $\mathrm{e}^{b_3}\le \beta_2$,~ $ \beta_2\mathrm{e}^{b_2}\le \beta_1 \mathrm{e}^{-b_1} $ and 
\begin{align}
    (\varepsilon^*)^2=\min &\left\{\frac{b_1}{6\mathrm{e}^{b_1}\breve K},~\frac{b_2\beta_2}{6\beta_1 h \breve L}, ~\frac{b_3 }{6\beta_1\eta\breve L}\right\}.
\end{align}
To maximize the value of $\varepsilon^*$, we choose  $\beta_1=\beta_2\mathrm{e}^{b_1+b_2}$ and $\beta_2=\mathrm{e}^{b_3}$, which yields
\begin{align}
    (\varepsilon^*)^2=\min &\left\{\frac{b_1}{6\mathrm{e}^{b_1}\breve K},~\frac{b_2 }{6\mathrm{e}^{b_1+b_2}h\breve L},~\frac{b_3}{6\mathrm{e}^{b_1+b_2+b_3}\eta \breve L}\right\}.
\end{align}
If we select $\varepsilon<\varepsilon^*$, there exists a $ \alpha_0(\varepsilon)>0$, such that 
\begin{align*}
    \dot V\leq -\alpha_0 V,
\end{align*} 
where 
\begin{align*}
   \alpha_0=\min  &\left\{b_1-6\varepsilon^2\mathrm{e}^{b_1}\breve K,~~
   \frac{b_2}{h}-6\varepsilon^2h\mathrm{e}^{b_1+b_2} \breve L, \right. \\  &  \left. ~\frac{b_3}{\eta}-6\varepsilon^2\eta\mathrm{e}^{b_1+b_2+b_3} \breve L\right\}
\end{align*}
which yields $V\le V(0)\mathrm{e}^{-\alpha_0t}$.
It is derived from \eqref{def-V12}-\eqref{def-V56},
 \begin{align}
   m_1(V_1+V_3+V_5)\leq V \leq m_2 (V_1+V_3+V_5),
 \end{align}
 with 
 \begin{align*}
   m_1=\min\{\beta_1 \mathrm{e}^{-b_1},\beta_2 h, \eta \},~~
   m_2=\max\{\beta_1, h\beta_2 \mathrm{e}^{b_2},~\eta  \mathrm{e}^{b_3}\}.
 \end{align*}
 Therefore
  \begin{align*}
     (V_1+V_3+V_5)\leq \frac{m_2}{m_1}(V_1(0)+V_3(0)+V_5(0))\mathrm{e}^{-\alpha_0t}.
 \end{align*}
Also, we get the $L^2$ norm relationship between the states of \eqref{eq:main-x-1}-\eqref{eq:bnd-u-1} and those  of \eqref{hat-target-z}-\eqref{eq:bnd-u1},
  		\begin{align*}
   \frac{1}{m_4}V_{0}\leq V_1+V_3+V_5\leq m_3 V_0,
  		\end{align*}
  		where $V_0=\|x\|^2_{L^2}+\|v\|^2_{L^2} +\|u\|^2_{L^2}$, with 
  		\begin{align*}
   m_3&=\max\{4(1+\bar K^2), 1+4h^2\bar L^2, 1+4\eta ^2 \bar J^2\},\\
   m_4&=\max\{4(1+\bar B^2), 1+4\bar D^2, 1+4 \bar E^2\}.
  		\end{align*}
Hence, we arrive at the stability bound \eqref{Th-feedback-ineq} with 
\begin{align*}
  W_0=\frac{m_2m_3m_4}{m_1}.
\end{align*} 
\end{proof}

\subsection{Stabilitzation of the observer error system under DeepONet observer gain}\label{NO-observer-stab}

\color{black}
\begin{theorem}\label{Th-Nob-sta}
   \rm{
For any  $( h, f) \in \mathbb{R}^+\times C^1(\mathcal{T}_1) $, there exist a sufficiently small $\varepsilon^*>0$, such that observer  
\begin{align}
   \hat x_t(s,t)=&~-\hat x_{s}(s,t)+\int_{s}^1  f(s,q)\hat x(q,t)dq
   \label{eq:main-x-oN}\\\nonumber
   &~+c(s)\hat u (0,t)+\hat Q_{1}(s)(v(0,t)-\hat v(0,t)),~\\ 
   \hat x(0,t)= &~U(t),\\
   \label{eq:main-v-oN}   h\hat v_t(s,t)= &~\hat v_s(s,t)
   +\hat Q_{2}(s)(v(0,t)-\hat v(0,t)),
   \\\label{eq:bnd-v-oN} \hat v(1,t)=&~\hat x(1,t),\\
   \label{eq:main-u-oN} 
   \eta \hat u_t(s,t)
   = &~\hat u_s(s,t),
   \\\label{eq:bnd-u-oN} \hat u(1,t)= &~x(1,t-h),
  		\end{align}
with all NO observer gains $\hat Q_i:=\hat{\mathcal{Q}}_i(h,f)$, defined in \eqref{ope-Qi}, $i=1,2$ of approximation accuracy $\varepsilon \in (0,\varepsilon^*)$ ensures that the observer error system, for all initial conditions $x_0,~\hat x_0,~ v_0, ~\hat v_0,~u_0,~\hat u_0 \in L^2[0,1]$, satisfies the exponential stability bound  
\begin{align}\label{error-st-bound}
   &~~\|x-\hat x\|^2_{L^2}+\|v-\hat v\|^2_{L^2}+\|u-\hat u\|^2_{L^2} \\\nonumber 
   &\leq W_1 \mathrm{e}^{-\alpha_1t} \left(\|x_0-\hat x_0\|^2_{L^2}+\|v_0-\hat v_0\|^2_{L^2}+\|u_0-\hat u_0\|^2_{L^2}\right), 
\end{align}
  		with $W_1>0$ and $\alpha_1>0$. 
}
\end{theorem}

\begin{proof}
Before proceeding, let $\tilde Q_i=Q_i-\hat Q_i$,  $i=1,2$ denote the difference between the exact observer gain and the neural operators. Similar to the proof of Theorem \ref{th-fb-st}, the proof contains two steps. 
  		First, we employ the transformation \eqref{eq:obs-transform-z-x} and \eqref{eq:obs-transform-w-v} to convert the error system \eqref{eq:main-x-e}-\eqref{eq:bnd-u2-e}, where the gains $Q_i$ are replaced with the NO observer gain $\hat Q_i$, $i=1,2$, to the following target system
  		\begin{align}
   \label{eq:target-zN}
   \tilde z_t(s,t)&=-\tilde z_s(s,t)+S(s)\tilde u(0,t)+\delta_1(s) \tilde w(0,t),\\
   \label{eq:target-bnd-zN} 
   \tilde  z(0,t)&=0,\\
   \label{eq:target-wN}
   h\tilde w_t(s,t)&=\tilde w_s(s,t)+\delta_2(s) \tilde w(0,t),
   \\\label{eq:target-bnd-wN}\tilde w(1,t)&=\tilde z(1,t),
   \\
   \eta \tilde u_t(s,t)&= \tilde u_s(s,t),
   \\\label{eq:target-bnd-u2-oN}\tilde u (1,t)&= 0,
\end{align}
 where 
\begin{align}
\delta_1(s,q)=&\int_s^1 F(s,q)\delta_1(q)dq +\frac{1}{h}\int_0^s M(s,q) \delta_2(q)dq,\nonumber \\
   &+\frac{1}{h}\int_s^1 P(s,q) \delta_2(q)dq+\tilde Q_1 (s)\label{delta4}
   \\
   \delta_2(s)=& \int_0^s R(s-q)\delta_2(q)dq+\tilde Q_2 (s),\label{delta5}
  		\end{align}
    and $S(s)$ is defined in \eqref{Eq:S}.
  		With \eqref{DON-Qi}  in Theorem  \ref{Th-NObs}, it is obvious that 
\begin{align}
   \|\delta_1\| \leq &\bar \delta_1:=\varepsilon\gamma_0,\label{bar-delta4}\\
   \|\delta_2\| \leq &\bar \delta_2:=\varepsilon\mathrm{e}^{\bar M}
   ,\label{bar-delta5}
\end{align}
    where $\gamma_0 =\mathrm{e}^{\bar F}\left(\frac{2\bar M \mathrm{e}^{\bar M}}{h}+1\right)$.
    
Second, we introduce the Lyapunov functional
  		\begin{align}
   V_{10}=\beta_3 V_7+V_8+\beta_4 V_9,
  		\end{align}
  		with 
  		\begin{align}
   V_7&=\int_0^1\mathrm{e}^{-b_4 s}\tilde z^2(s,t) ds,\\
   V_8&=  h\int_0^1\mathrm{e}^{b_5 s}\tilde w^2(s,t) dr,~~
   V_9=\eta \int_0^1  \mathrm{e}^{b_6 s}\tilde u^2(s,t) dr,
  		\end{align}
  		and $\beta_3,~\beta_4$ are positive constants.
  		Taking the time derivative, we get
  		\begin{align}
   \dot V_{10}=&~\beta_3 \int_0^1 -\mathrm{e}^{-b_4 s}(\tilde z^2(s))_sds+\int_0^1 \mathrm{e}^{b_5 s}(\tilde w^2(s))_s ds\nonumber \\&
   +2\beta_3 \int_0^1 \mathrm{e}^{-b_4s}\tilde z(s)S(s)ds \tilde u(0)
   \nonumber \\
   &+2\beta_3  \int_0^1 \mathrm{e}^{-b_4s}\tilde z(s)\delta_1(s)ds \tilde w(0)
   \nonumber \\
   &+2  \int_0^1 \mathrm{e}^{b_5s}\tilde w(s,t)\delta_2(s)ds \tilde w(0)
   \nonumber\\&+\beta_4 \int_0^1 \mathrm{e}^{b_6 s}(\tilde u^2(s))_s ds,
  		\end{align}
  		where we have used $\tilde v(0)=\tilde w(0)$ from \eqref{eq:obs-transform-w-v}. 
  		\begin{align}
   \dot V_{10}\leq &-\beta_3 \mathrm{e}^{-b_4}\tilde z^2(1)-\beta_3 b_4 V_7 +
   \mathrm{e}^{b_5}\tilde w(1)-\tilde w^2(0)-\frac{b_5}{h} V_8\nonumber \\
   &-\beta_4 \tilde u^2(0)-\frac{1}{\eta}\beta_4b_6   V_9+\beta_3 \bar S V_7+\beta_3 \bar S \tilde u^2(0)\nonumber \\
   &+\beta_3 \bar \delta_1 V_7 +\beta_3 \bar \delta_1 \tilde w^2(0)+\frac{\bar \delta_2}{h} V_8+\delta_2 \tilde w^2(0)\nonumber \\ 
   \leq &-(\beta_3 \mathrm{e}^{-b_4}-\mathrm{e}^{b_5})\tilde z^2(1)-(1-\beta_3\bar \delta_1-\bar \delta_2)\tilde w^2(0)
   \nonumber \\&
   -(\beta_4-\beta_3\bar S)\tilde u^2(0)-\beta_3( b_4 -\bar S -\bar \delta_1)V_7\nonumber \\& -\frac{1}{h}(b_5-\bar\delta_2 ) V_8
   -\frac{1}{\eta } \beta_4b_6 V_9. 
\end{align}
We choose $\beta_3 =\mathrm{e}^{b_4+b_5}$, $\beta_4\ge \beta_3\bar S=\bar S\mathrm{e}^{b_4+b_5}$ and $b_4>\bar S$, such that $\varepsilon^*$ as large as possible,  
\begin{align}
    \varepsilon^* =\min &\left\{\frac{1}{ \mathrm{e}^{b_4+b_5}\gamma_0+\mathrm{e}^{\bar M}},
    \frac{b_5}{\mathrm{e}^{\bar M}},
    \frac{b_4-\bar S}{\gamma_0} \right\}.
\end{align}
If appropriate parameters $\beta_4$, $b_4$ and $b_5$ are selected, such as $b_4=2\bar S$, and $\beta _4=2\bar S\mathrm{e}^{2\bar S+b_5}$, we get
\begin{align}
    \dot V \leq & -(1-\varepsilon(\beta_3\gamma_0+\mathrm{e}^{\bar M}))\tilde w^2(0)-\bar S\mathrm{e}^{2\bar S+b_5}\tilde u^2(0)-\frac{1}{\eta } \beta_4b_6 V_9
   \nonumber \\&
   -\beta_3( \bar S -\varepsilon\gamma_0)V_7  -\frac{1}{h}(b_5-\varepsilon\mathrm{e}^{\bar M} ) V_8. 
\end{align}
Therefore, there exists a $ \alpha_1(\varepsilon)>0$ for $\varepsilon<\varepsilon^*$, such that 
\begin{align*}
    \dot V_{10}\leq -\alpha_1 V_{10},
\end{align*}    
with 
\begin{align}
   \alpha_1=\min\left\{ b_4-\bar S-\varepsilon \gamma_0 ,\frac{1}{h}(b_5-\varepsilon \mathrm{e}^{\bar M}), \frac{1}{\eta}  b_6\right\},
\end{align}
\color{black}
namely, $V_{10}\leq V_{10}(0) \mathrm{e}^{-\alpha_1t}$.  Denoting $V_{11}:=\|\tilde z \|^2_{L^2}+\|\tilde v \|^2_{L^2}+\|\tilde u \|^2_{L^2}$, we have   
\begin{align}
   m_5V_{11} \leq V_{10} \leq m_6 V_{11},
  		\end{align}
  		with 
  		\begin{align*}
   m_5&=\min\{\beta_3 \mathrm{e}^{-b_4}, ~h, ~ \eta\beta_4 \},\\
   m_6&=\max\{\beta_3, ~h \mathrm{e}^{b_5},~\eta  \beta_4\mathrm{e}^{b_6}\}.
\end{align*}
Given that transformation \eqref{eq:obs-transform-z-x}-\eqref{eq:obs-transform-w-v} is invertible, with the inverse transformation  defined in \eqref{eq:obs-transform-x-z}-\eqref{eq:obs-transform-v-w}, and the bounds of the kernels are presented in Theorem \ref{Th-obskernel}, norm equivalence holds between the observer error system and the associated target system in the following sense:
\begin{align*}
   \frac{1}{m_7}V_{11}\leq \|\tilde  x\|^2_{L^2}+\|\tilde v\|^2_{L^2}+\|\tilde u\|^2_{L^2} \leq m_{8} V_{11},
  		\end{align*}
  		with 
  		\begin{align*}
   m_7&=\max\{4+4\bar F^2, 2+10h^2\bar F^2, 1\},\\
   m_{8}&=\max\{4+4\bar F^2, 2+10h^2\bar M^2, 1\}.
  		\end{align*}
Hence, we arrive the stability bound \eqref{error-st-bound} with 
\begin{equation*}
 W_1=\frac{m_6m_7m_8 }{m_5}.
\end{equation*}
\end{proof}

\subsection{Output-feedback stabilization with DeepONet gains for controller and observer} 
\label{output-feedback} 
\begin{figure}[t]    
  		\centering
  		\includegraphics[width=0.45\textwidth]{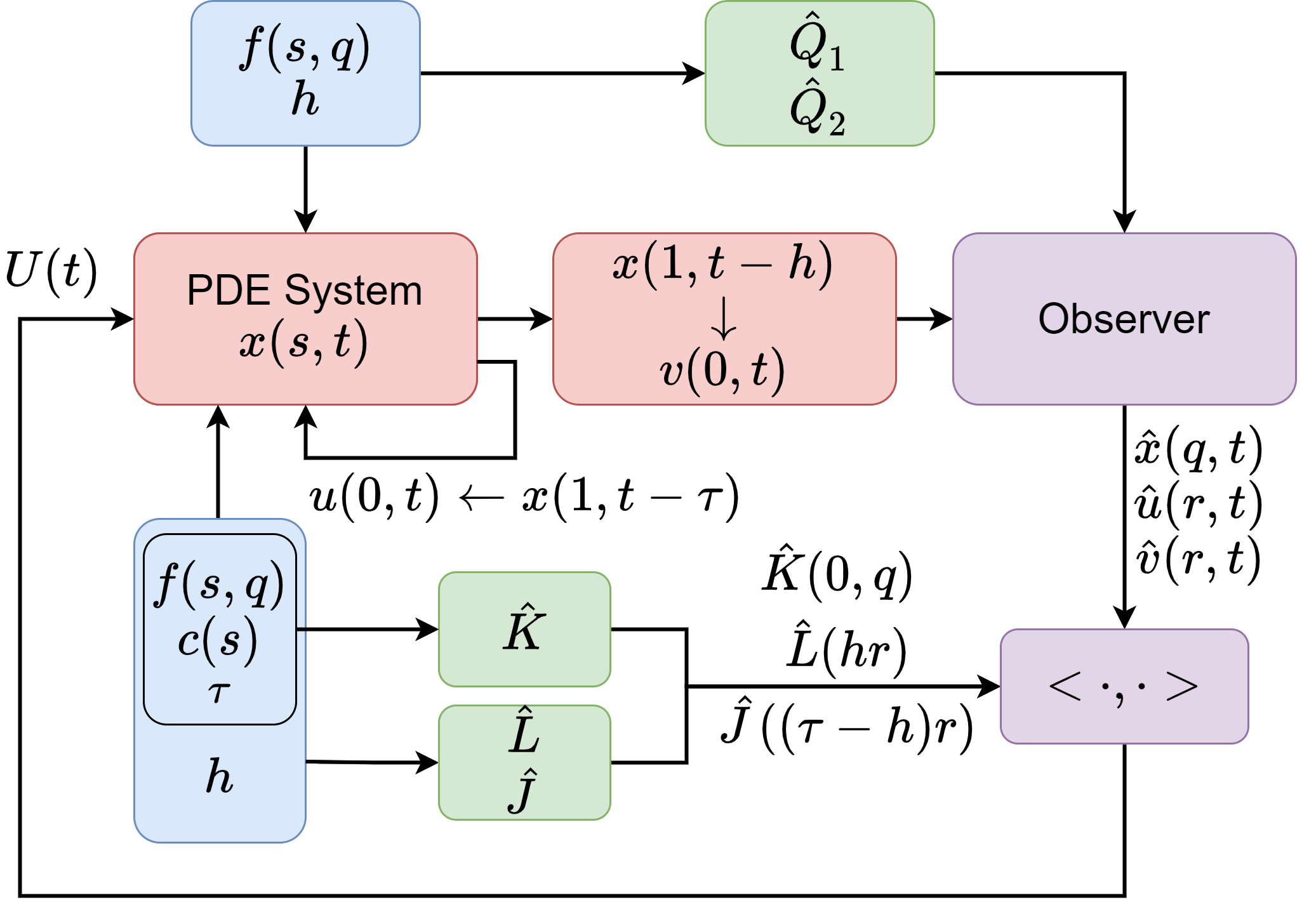}
  		\caption{The neural operator learning framework for backstepping delay compensation control.}
  		\label{fig:system}
\end{figure}

   In this section, we put together the observer \eqref{eq:main-x-oN}-\eqref{eq:bnd-u-oN}  along with the observer-based controller 
  	\begin{align}\label{estimated-controll}
  		U(t)=&\int_0^1  \hat K(0,q) \hat x(q,t)dq+h\int_0^1   \hat L(hr)\hat v(r,t)dr
  		\nonumber \\&+\eta \int_0^1 \hat J(\eta r)\hat u(r,t)dr ,
  	\end{align}
  	to stabilize the system \eqref{eq:main-x-1}-\eqref{eq:bnd-u-1}. 
  	Figure \ref{fig:system} illustrates the framework of the neural operator based output feedback for the delayed PDE system. As shown in Figure \ref{fig:system}, We apply three neural operators to learn the operators $K$, $L$ and $J$ defined in \eqref{ope-K}-\eqref{ope-J}, then to derive the gain functions which are used in the controller. For the observer, we apply two neural operators to learn the operator $Q_1$ and $Q_2$ defined in \eqref{ope-Qi}, which are used in the observer. We use the estimated system states for feedback with the learned neural gain functions in the control law. The control kernel  and the observer gain functions can be learned once. The trained DeepONets are ready to produce the control kernel and observer gain functions for any new functional coefficients and any new delays.   
  	
The following theorem establishes the exponentially stability for the cascading system under the output-feedback control with the DeepONet gains. 
  	\begin{theorem}\label{th-output-st}
   \rm{
  		Consider the system \eqref{eq:main-x-1}-\eqref{eq:bnd-u-1}, along with the observer \eqref{eq:main-x-oN}-\eqref{eq:bnd-u-oN} and the control \eqref{estimated-controll}, where the exact backstepping control kernels $ K$, $ L$, $ J$ and observer gains $ Q_1$, $ Q_2$ are approximated by DeepONets $\hat K$, $\hat L$, $\hat J$ and $\hat Q_1$, $\hat Q_2$, respectively with the accuracy $\varepsilon\in (0,\varepsilon^*)$.
  		For any  $(\tau, f, c) \in \mathbb{R}^+ \times C^1(\mathcal{T}_1) \times \underline C $ corresponding to the control kernels $ K$, $ J$, $(\tau, \eta, f, c) \in \mathbb{R}^+ \times \mathbb{R}^+ 
 \times C^1(\mathcal{T}_1) \times \underline C $, to $L$,  and $( h, f) \in \mathbb{R}^+\times C^1(\mathcal{T}_1)$, corresponding to observer gains $Q_1$, $ Q_2$, there exist a sufficiently small $\varepsilon^*>0$, such that the observer-based control \eqref{estimated-controll} 
  		ensures that the observer cascading closed-loop system satisfies the exponential stability bound, for all $t>0$
  		\begin{align}\label{Th-output-ineq}
   \Theta(t) \leq  W_2 \mathrm{e}^{-\alpha_2 t} \Theta(0),  
  		\end{align}
  		where 
  		\begin{align*}
   \Theta(t)=\| x\|^2_{L^2}+\| v\|^2_{L^2} +\| u\|^2_{L^2}+\|\hat x\|^2_{L^2}+\|\hat v\|^2_{L^2} +\|\hat u\|^2_{L^2},
  		\end{align*}
  		with $W_2>0$ and $\alpha_2>0$. 
}
\end{theorem}

\begin{proof}
We consider the observer error system
    \begin{align}
   \tilde x_t(s,t)=&~-\tilde x_{s}(s,t)+
   \int_{s}^1  f(s,q)\tilde x(q,t)dq\nonumber
   \\ &~ +c(s)\tilde u(0,t)-\hat Q_{1}(s)\tilde v(0,t),\label{eq:main-x-eN}\\ \label{eq:bnd-x-eN}
   \tilde x(0,t)=&~0,\\\label{eq:main-u-eN}
   h \tilde  v_t(s,t)=&~  \tilde v_s(s,t)-\hat Q_{2}(s)\tilde v(0,t),
   \\\label{eq:bnd-u-eN}\tilde v(1,t)=&~\tilde x(1,t),
   \\
   \eta \tilde u_t (s,t)= &~\tilde u_s (s,t),
   \\\label{eq:bnd-u2-eN}  \tilde u(1,t)=&~ 0, 
\end{align} 
and the observer \eqref{eq:main-x-oN}-\eqref{eq:bnd-u-oN} with control \eqref{estimated-controll}, since they are equivalent to the  cascading system  \eqref{eq:main-x-1}-\eqref{eq:bnd-u-1}  and observer \eqref{eq:main-x-oN}-\eqref{eq:bnd-u-oN} with control \eqref{estimated-controll}.
  		
The proof contains two steps.
First, we derive the target system of $(\hat x, \hat v, \hat u)-(\tilde x, \tilde v, \tilde u)$ by applying the backstepping transformation 
  		\begin{align}
   \breve z(s,t)= &~\Gamma|_{  K,  L, J}[\hat x,\hat v,\hat u](s,t),   \label{eq:trans-hat-breve}\\
   \tilde x(s,t) = &~\mathcal{F}|_{F,M,P}[\tilde z, \tilde w](s,t)   \label{eq:trans-tilde-zx},\\
   \tilde v(s,t) = &~\mathcal{R}|_{R}[\tilde w](s,t) , \label{eq:trans-tilde-wv}
  		\end{align}
combining the DeepONet approximated control gains $\hat K(0,t)$, $\hat L(hr)$ and $\hat J(\eta r)$, and observer gains $\hat Q_1(s)$ and $\hat Q_2(s)$ to transform the observer \eqref{eq:main-x-oN}-\eqref{eq:bnd-u-oN} cascading observer error system 
  		into the $(\breve z, \hat v, \hat u)-(\tilde z, \tilde w, \tilde u)$ system as 
  		\begin{align}
   \breve z_t(s,t)=&~- \breve z_{s}(s,t)+ G(s)\tilde w(0),\label{breve-target-z}\\ 
   \label{breve-bnd-z}
   \breve z(0,t)= &~-\int_0^1 \breve z(q)[\tilde K(0,q)+\int_0^q \tilde K(0,r)B(r,q)dr]dq\nonumber 
    \\&-\int_0^1 \hat v(q)[h\tilde L(hq)+\int_0^1 \tilde K(0,r)D(r,q)dr]dq\nonumber 
    \\&-\int_0^1 \hat u(q)[\eta\tilde J(\eta q)+\int_0^1 \tilde K(0,r)E(r,q)dr]dq,\\ \label{eq:tar-v}
   h \hat v_t(r,t)=&~\hat v_s(r,t)+\hat Q_2(s)\tilde w(0,t), 
   \\\label{eq:tar-bnd-v}
   \hat v(1,t)=&~\breve z(1,t),\\
   \label{eq:tar-main-u}
   \eta \hat u_t(r,t)=&~\hat u_s(r,t),\\\label{eq:tar-bnd-u}
   \hat u(1,t)=& ~v(0,t)=\hat v(0,t)+\tilde w(0,t),
\end{align}
where $\delta_i$, $i=1,2$ are given in \eqref{delta4}-\eqref{delta5} and
\begin{align}
     G(s)=& ~\hat Q_1(s)-\int_s^1   K(s,q)\hat Q_1(q)dq\nonumber \\&- \int_0^1  L(s+hr)\hat Q_2(r)dr.
\end{align}
\color{black}
Since kernels $\hat Q_1$, $\hat Q_2$, $K$ and $L$ are given bounded,  
in Theorem \ref{Th-kernel}, \ref{Th-obskernel} and \ref{Th-NObs}, 
 $G(s)$ is also bounded, denoting the bound of $|G|$ by $\bar G$.
It worth noting that we still employ the transformations involving the analyzed control and observer kernel functions, albeit with the substitution of control gains and observer gains by the DeepONets in the control  and observer implementation.    

Second, we introduce the following Lyapunov functional to prove the stability of the cascading target system \eqref{breve-target-z}-\eqref{eq:tar-bnd-u} and \eqref{eq:target-zN}-\eqref{eq:target-bnd-u2-oN}.
Since target system \eqref{breve-target-z}-\eqref{eq:tar-bnd-u} has the same form as that of the target system of the state-feedback system, except for one extra term in \eqref{breve-target-z} and in \eqref{eq:tar-v}, respectively, and boundary condition \eqref{eq:tar-bnd-u}, we redefine the Lyapunov functionals \eqref{def-V12}-\eqref{def-V56} as follows:
  		\begin{align} \label{def-V12O}
   V_1&=\|   \breve z^2\|^2_{L^2}, ~~
   V_2 =\int_0^1 \mathrm{e}^{-b_1 s}\breve z^2(s,t)ds, \\
   \label{def-V34O}
   V_3&=\|   \hat v^2\|^2_{L^2}, ~~
   V_4=h\int_0^1 \mathrm{e}^{b_2s}  \hat v^2(s,t)ds,\\
   \label{def-V56O}
   V_5&= \|   \hat u^2\|^2_{L^2}, ~~
   V_6=\eta \int_0^1 \mathrm{e}^{b_3s}\hat u^2(s,t)ds,
  		\end{align}
  		with $b_i>0$, $i=1,2,3$.
Take time derivative of $V:=\beta_1 V_2+\beta_2 V_4 +V_6+\beta_5 V_{10}$ with $\beta_i>0$, $i=1,2,5$, we have 
\begin{align}
   \dot V=& -\beta_1 \mathrm{e}^{-b_1} \breve z^2(1)+\beta_1 \breve z^2(0)-\beta_1 b_1 V_2 +\beta_2  \mathrm{e}^{b_2} \hat v^2(1)\nonumber\\&
     -\beta_2 \hat v^2(0) +2\beta_1 \int_0^1 \mathrm{e}^{-b_1 s}\breve{z}(s)G(s)ds \tilde{w}(0)\nonumber \\& + \mathrm{e}^{b_3} (\hat v(0)+\tilde w(0))^2+2\beta_2 \int_0^1 \mathrm{e}^{b_2 s}\hat{v}(s)\hat Q_2(s)ds \tilde{w}(0) \nonumber \\
   &-\beta_2 b_2  V_4/h-\hat u^2(0)  
    -b_3 V_6/\eta+\beta_5 \dot V_{10}.
   \label{R-derivative-V}
  		\end{align}
and combining \eqref{ieq:z0}, we get 
\begin{align}
   \dot V\le & -(\beta_1\mathrm{e}^{-b_1}-\beta_2\mathrm{e}^{b_2}) \breve z^2(1)
   -(\beta_2- 2\mathrm{e}^{b_3}) \hat v^2(0) -\hat u^2(0)
   \nonumber \\&
   -\beta_5(\beta_4-\beta_3 \bar S) \tilde u^2(0)-\beta_5(\beta_3 \mathrm{e}^{-b_4}-\mathrm{e}^{b_5})\tilde z^2(1)\nonumber \\&
   -\left(\beta_5-2\mathrm{e}^{b_3}-\beta_1\bar G-\beta_2\mathrm{e}^{b_2}\bar M-\varepsilon(\beta_2\mathrm{e}^{b_2}+\beta_3\beta_5\gamma_0 \right.
   \nonumber \\
   &\left.+\beta_5 \mathrm{e}^{\bar M})\right)\tilde w^2(0)-\beta_1(b_1-\bar G-6\varepsilon^2\breve K\mathrm{e}^{b_1})V_2\nonumber 
   \\&-(\beta_2(b_2/h-\bar M-\varepsilon)-6\varepsilon^2\beta_1h\breve L)V_4\nonumber \\&-(b_3/\eta-6\varepsilon^2\beta_1\eta \breve L)V_6 -\beta_5\beta_3( b_4 -\bar S -\varepsilon\gamma_0 )V_7\nonumber \\& -\frac{\beta_5}{h}(b_5-\varepsilon\mathrm{e}^{\bar M}  ) V_8
   -\frac{\beta_5}{\eta } \beta_4b_6 V_9   ,    
\end{align}
where $\breve K$ and 
$\breve L$ are defined in \eqref{define-breve_K_L}. To maximize the approximation error
$\varepsilon^*$, we choose 
\begin{align}
    \beta_3=\mathrm{e}^{b_4+b_5},~~\beta_4=\mathrm{e}^{b_4+b_5}\bar S,~~\beta_1=\beta_2 \mathrm{e}^{b_1+b_2},
\end{align}
and based on the above values, one can further determine the value of $\beta_2=2\mathrm{e}^{b_3}$, which gives
\begin{align}
    \varepsilon^*=&\min\left\{\frac{\beta_5-2\mathrm{e}^{b_3}(1+\mathrm{e}^{b_1+b_2}\bar G+\mathrm{e}^{b_2}\bar M)
    }{2\mathrm{e}^{b_2+b_3}+\beta_5(\gamma_0 \mathrm{e}^{b_4+b_5}+\mathrm{e}^{\bar M})},\right.\nonumber \\
    & \left.\frac{\sqrt{b_1-\bar G}}{\sqrt{6\mathrm{e}^{b_1}\breve K}},~\frac{b_3}{12\mathrm{e}^{b_1+b_2+b_3}\eta^2\breve L},~\frac{b_4-\bar S}{\gamma_0},\right.
    \nonumber \\&\left.\frac{\sqrt{1+24\mathrm{e}^{b_1+b_2}h\breve L(b_2/h-\bar M)}-1}{12\mathrm{e}^{b_1+b_2}h\breve L},~\frac{b_5}{\mathrm{e}^{\bar M}}
    \right. \Bigg\}.
\end{align}
Selecting appropriate value of $\beta_5$, $b_1$, $b_2$ and $b_4$, one can get  $0<\varepsilon<\varepsilon^*$ such that $\dot V \leq -\alpha_2(\varepsilon) V$,
with 
\begin{align}
    \alpha =\min&\left\{ b_1-\bar G-6\breve K\mathrm{e}^{b_1}\varepsilon,~b_2/h-\bar M-\varepsilon-6
    \mathrm{e}^{b_1+b_2}h\breve L\varepsilon^2,  \right.
\nonumber \\& \left.~b_3/\eta-12\mathrm{e}^{b_1+b_2+b_3}\eta \breve L\varepsilon^2,
 b_4-\bar S-\gamma_0\varepsilon,\right.
\nonumber \\& \left.
~1/h(b_5-\mathrm{e}^{\bar M}\varepsilon),  1/\eta  \bar Sb_6 \right\}.
\end{align}
\color{black}
Hence $V\leq V(0)\mathrm{e}^{-\alpha_2 t}$.
Due to $L^2$ norm equivalence which are proven in Theorem \ref{th-fb-st}  and \ref{Th-Nob-sta}, we conclude that there exist a positive constant $W_2$ such that 
\eqref{Th-output-ineq} holds.
\end{proof}
  	
  	Throughout the stability analysis of the overall cascading system, the "separation principle" can also be proven.

  	\section{Numerical Results: Full-State Feedback, Observer, and Output Feedback}\label{simulation}
  	
\subsection{Full-state feedback}
  	We first solve equations \eqref{eq:K}-\eqref{eq:J} numerically by using the finite difference method to get datasets for different delays $\tau$, $ h$ and different coefficient functions $f(s,q)$, $c(s)$ to train neural operators $\hat{\mathcal{K}}( \tau, f,c)$, $\hat{\mathcal{L}}( \tau,h, f,c)$ and $\hat{\mathcal{J}}( \tau,h, f,c)$. Let state delay $ \tau \sim U (0.8, 2)$,  sensor delay  $h \sim U(0.1,0.7)$, function $f(s,q)$ as  a product of Chebyshev polynomials $f(s,q) =9\cos\big(\mu_1 \cos^{-1}(s)\big)\cos\big(\mu_2 \cos^{-1}(q)\big)$ with $\mu_1, \mu_2 \sim U (3, 6)$, and function $c(s)$ as a Chebyshev polynomial $c(s)=\cos\big(\mu_3 \cos^{-1}(s)\big)-\cos\big(\mu_3 \cos^{-1}(1)\big)$ with $\mu_3\sim U (3, 6)$, where  $U (a, b)$ denotes the uniform distribution in the interval $[a,~b]$. In sampling, the discretized spatial step size is set to $\Delta s=0.02$. The Simulation code is shared on \href{https://github.com/JingZhang-JZ/NO_hyperbolic_delay.git}{github}.

  \begin{figure*}[!htbp]    
   \centering
   \includegraphics[width=0.7\textwidth]{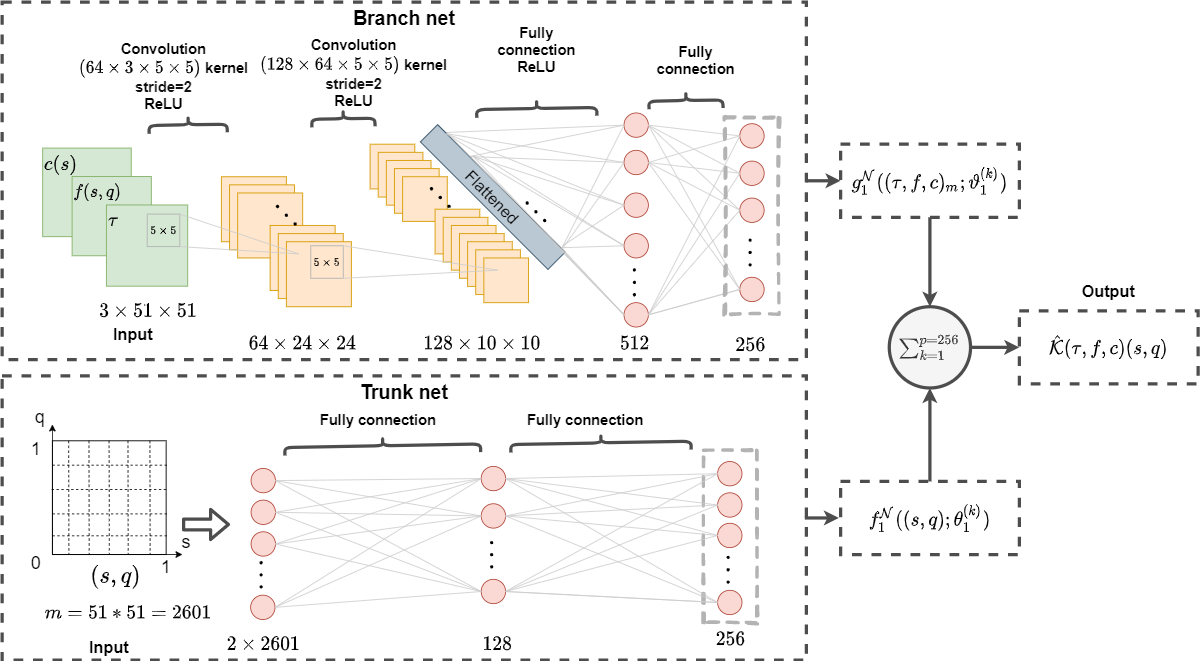}
   \caption{The DeepONet structure for kernel $K$.}
   \label{fig:branch}
\end{figure*}
 As shown in Figure \ref{fig:branch}, we construct a branch net consisting of  two layers  $5\times 5$ convolutional neural networks (CNNs) with  strides of $2$,   and  two layers of  $12800\times512$ and $512\times256$ fully connected networks,   and a trunk net consisting of two layers of $2601\times 128$ and $128\times 256$  fully connected networks. Being different from the neural network of $K$, the input to the neural networks of $L$ and $J$ is $4\times 51\times 51$. 
  	Three DeepONets  are employed to learn the three kernel functions $K$, $L$ and $J$,  which contain 6928641, 6930241 and 6930241 parameters, respectively. 
 The loss function is chosen as the smooth $L_1$  \cite{girshick2015fast}, by using $\hat \rho \in \mathbb{R} $  to denote the prediction  and $\rho \in \mathbb{R}$ to denote the true value:
  		\begin{align}\label{SL1Loss}
   \mathrm{Loss}(\hat \rho-\rho)= 
   \left\{\begin{aligned} 
   	& 0.5(\hat \rho-\rho)^2, &\mathrm{if} ~\vert \hat \rho-\rho \vert < 1, \\
   	&\vert  \hat \rho-\rho \vert -0.5, &\mathrm{otherwise.}
   \end{aligned}\right.
  		\end{align}
    The smooth $L_1$ loss can be seen as exactly $L_1$ loss, but with the $|\hat \rho-\rho|$ portion replaced with a quadratic function such that its slope is $1$ at $|\hat \rho-\rho|=1$. The quadratic segment when $|\hat \rho-\rho|<1$ smooths the $L_1$ loss near $|\hat \rho-\rho|=0$, avoiding sharp changes in slope. The smooth $L_1$ combines the advantages of the $L_1$ and $L_2$ loss functions. When the difference between the prediction and the true value is large, the gradient value won't be too large; when the difference between the prediction and the true value is small, the gradient value is also small.        
    We also train the DeepONets using the $L_2$ loss function and find that using the smooth $L_1$ loss function exhibits better performance compared to the $L_2$ loss function.
    The evolution of loss over time of using both loss functions is shown in Figure 3, which illustrates that using the smooth $L_1$ loss has smaller errors and less fluctuations. Therefore, we will apply the smooth $L_1$ to train the networks. If not specifically pointed out, all the NNs in the following simulations are trained using the smooth $L_1$  loss function.
 
\begin{figure*} 
\centering
\begin{tabular}{ccc}
\includegraphics[width=0.3\textwidth]{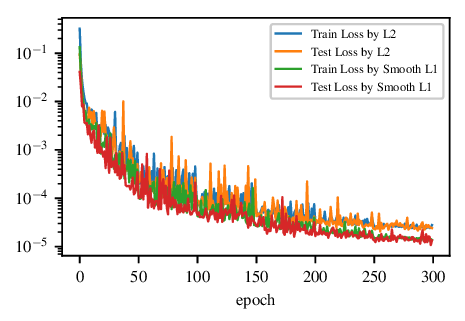}
&\includegraphics[width=0.3\textwidth]{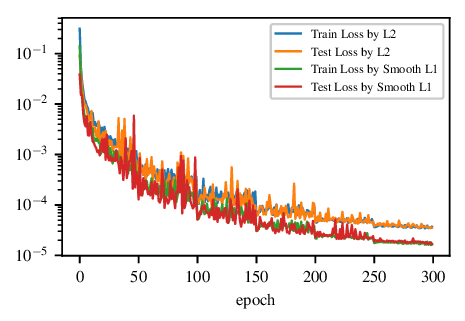}
&\includegraphics[width=0.3\textwidth]{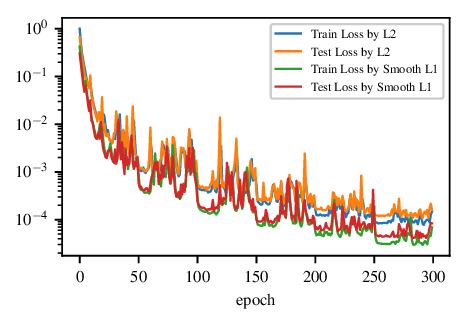}\\
(a) $\mathrm{Loss}(\hat K, K)$ 
& (b) $\mathrm{Loss}(\hat L, L)  + \mathrm{Loss}(\hat J, J)$ 
& (c) $\mathrm{Loss}(\hat Q_1, Q_1)  +  \mathrm{Loss}(\hat Q_2, Q_2)$ 
\end{tabular}
\caption{(a) The loss  of the neural control kernel for $K$. (b) The loss of the neural control kernels for $ L$ and $ J$. (c) The loss  of neural  observer gains.}
\label{fig:loss}
\end{figure*}
 	  	
\begin{figure*}
\centering
\begin{tabular}{ccc}
\includegraphics[width=0.31\textwidth,trim=0 150 162 5,clip]{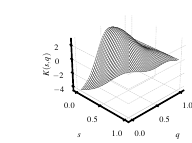}
&\includegraphics[width=0.31\textwidth,trim=0  150 162 5,clip]{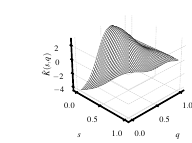}
&\includegraphics[width=0.31\textwidth,trim=0  150 162 5,clip]{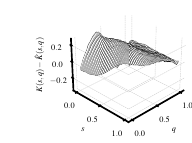}
\end{tabular}
\begin{tabular}{cc}
\includegraphics[width=0.41\textwidth]{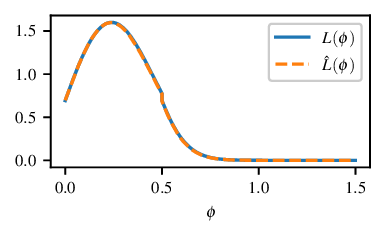}
&\includegraphics[width=0.41\textwidth]{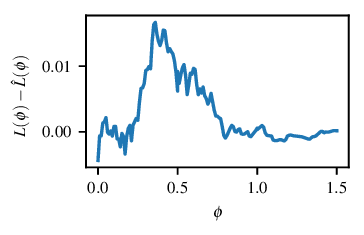} \\
\includegraphics[width=0.41\textwidth ]{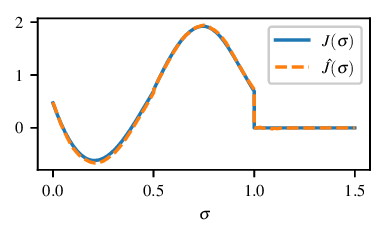}
&\includegraphics[width=0.41\textwidth]{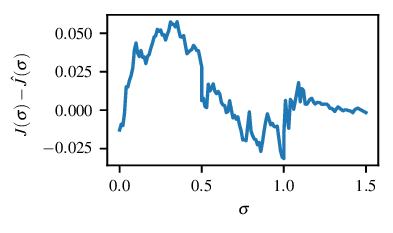}
\end{tabular}
\caption{The first row shows the kernel functions $K(s,q)$, the learned kernel functions $\hat K(s,q)$ and the errors $K(s,q)-\hat K(s,q)$.  The second row shows the kernel functions $L(s)$, the learned kernel functions $\hat L(\phi)$ and the errors $L(\phi)-\hat L(\phi)$. The last row shows the kernel functions $J(\sigma)$, the learned kernel functions $\hat J(\sigma)$ and the errors $J(\sigma)-\hat J(\sigma)$. }
\label{fig:K}
\end{figure*}

\begin{figure*}  
  		\centering
  		\begin{tabular}{ccc}
   \includegraphics[width=0.31\textwidth]{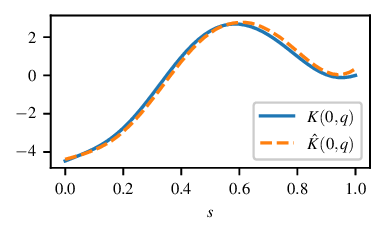}&
   \includegraphics[width=0.31\textwidth]{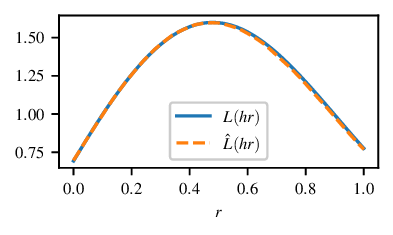}&
   \includegraphics[width=0.31\textwidth]{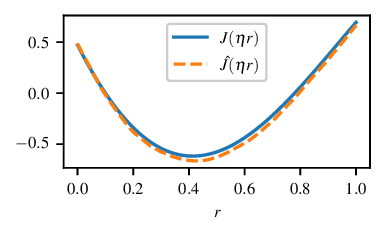}\\
   \includegraphics[width=0.31\textwidth]{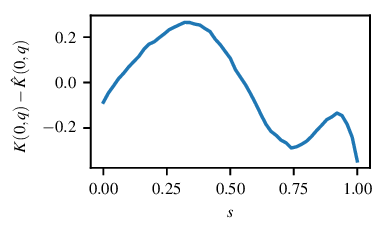}&
   \includegraphics[width=0.31\textwidth]{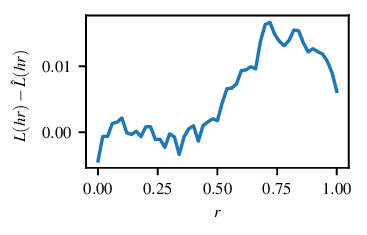}&
   \includegraphics[width=0.31\textwidth]{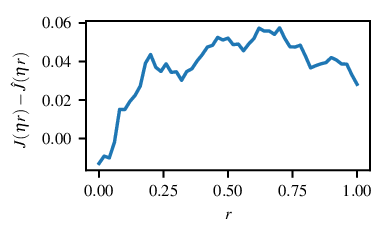}\\
  		\end{tabular}
  		\caption{The first row shows the analyzed   control gains $K(0,q)$, $L(hr)$, $J(\eta r)$, and the learned control gains $\hat K(0,q)$, $\hat L(hr)$, $\hat J(hr)$. The last row shows the errors $K_1(0,q)- \hat K_1(0,q)$, $L(hr)- \hat L(hr)$, $J(\eta r)- \hat J(\eta r)$.}
  		\label{fig:kbud}
  	\end{figure*}

Since operator $\mathcal{K}$'s input variables are three and the other operators have four input variables, we first train $\hat{\mathcal{K}}(\tau, f,c)$ on a dataset of  $8000$ numerical solutions with different parameters, and then train simultaneously the networks for operators $\hat{\mathcal{L}}( \tau,h, f,c)$ and $\hat{\mathcal{J}}( \tau,h, f,c)$  using  $10000$ instances. The NN for operator $K$ achieves a training loss of $1.36\rm{E}-5$ and a testing loss of $1.34\rm{E}-5$ after $300$ epochs in around $11$ minutes, shown in Figure~\ref{fig:loss} (a).  The  two NNs for operator $L$ and $J$ achieve  a training loss of $1.90\rm{E}-5$ and a testing loss of $1.90\rm{E}-5$ after $300$ epochs in  about $15$ minutes, which are shown in Figure~\ref{fig:loss} (b). (The experimental code runs on  Intel$^\circledR$ Core$^\text{TM}$ i9-7900X CPU @ 3.30GHz $\times$ 20 and GPU TITAN Xp/PCle/SSE2.)  	
  	In Figure~\ref{fig:K}, we demonstrate the analytical kernels that are solved numerically,  the learned DeepONet kernels and the errors between them, where the coefficients are chosen as  $ \tau=1$, $h=0.5$,  $f(s,q)$ with $\mu_1 =\mu_2 = 5$ and  $c(s)$  with $\mu_3=5$. Also,  in Figure~\ref{fig:kbud}, we show both the analytical control gains and learned control gains, respectively. 

\begin{figure}
\centering
\includegraphics[width=0.45\textwidth]{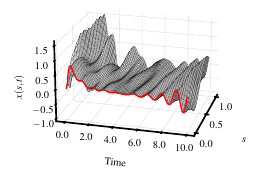}
\caption{The system dynamics  without delay compensation.}
\label{fig:uncom}
\end{figure}

  	To test the performance of the neural operator based control, we apply the trained neural gains in controller \eqref{hat-controller}. Here, we use the same parameter settings as Figure \ref{fig:K} and let initial condition be $x(s,0)=\sin(x)$. The upwind scheme with a time step size of $\Delta t=0.001$ and a trapezoidal integration rule are used to  numerically solve the PIDE under controller \eqref{hat-controller}. 
  	Before proceeding, we show in  Figure~\ref{fig:uncom} that the dynamical state $x(s, t)$ of the system under a nominal controller without delay compensation fails to converge.
  	
  	
\begin{figure*}
   
  		\centering
  		\includegraphics[width=0.95\textwidth,trim=0 130 0 30,clip]{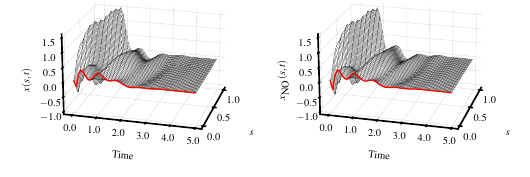}
  		\begin{tabular}{cc}
   \includegraphics[width=0.4\textwidth]{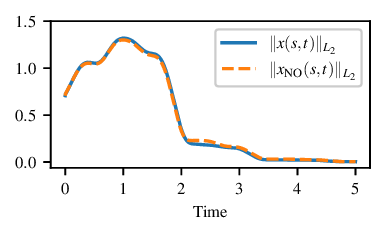}
   &\includegraphics[width=0.45\textwidth]{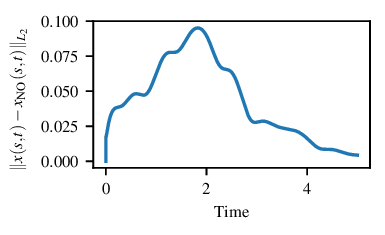}
  		\end{tabular}
  		\caption{ The closed-loop evolution with full-state feedback \eqref{hat-controller}.  The left column in the first row shows state $x(s,t)$ with the analyzed kernels $K$, $L$ and $J$. The right column in the first row shows state $x_{\mathrm{NO}}(s,t)$ with the learned kernels $\hat K$, $\hat L$ and  $\hat J$. The last row shows the $L_2$-norm of state $x$, $ x_{\mathrm{NO}}$, and the error between them.}
  		\label{fig:statefeedback}
  	\end{figure*}

In Figure~\ref{fig:statefeedback}, we demonstrate the dynamics of the closed-loop with the full state feedback, using the numerically solved control gains and the DeepONet learned control gains, respectively. The closed-loop system dynamics with NO kernels approximates the PDE well with a peak error of less than $8\%$  compared to the closed-loop system with analytical kernels. 
	
  	\subsection{Output feedback}
  	We train two neural observer gains $\hat{Q}_{1}(s)$ and $\hat{Q}_{2}(s)$ instead of the four observer kernels, which reduces the computational cost in half. 
  	The same parameter settings as for the full-state feedback are applied in the NN training, and the sensor delay $h$ is chosen from $U (0.1, 0.6) $.  
  	
\begin{figure*}    
  		\centering
  		\begin{tabular}{cc}
   \includegraphics[width=0.4\textwidth]{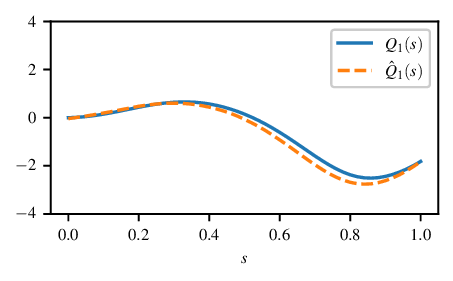}&
   \includegraphics[width=0.4\textwidth]{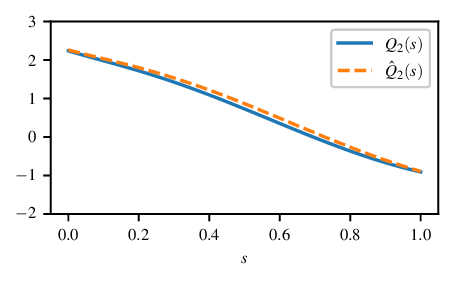}\\
   \includegraphics[width=0.4\textwidth]{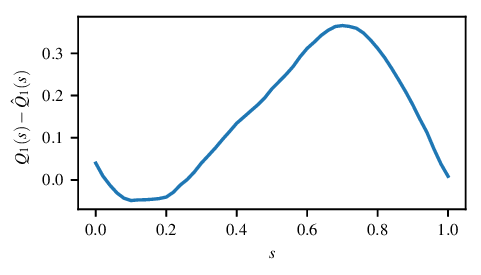}&
   \includegraphics[width=0.4\textwidth]{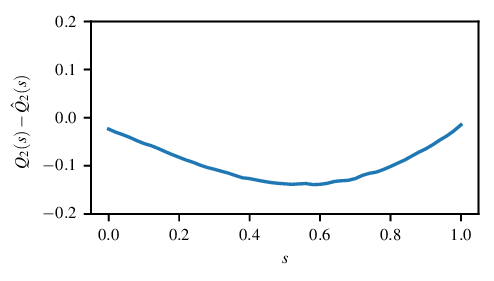}
  		\end{tabular}
  		\caption{The first row shows the analyzed observer gains $Q_1(s)$, $Q_2(s)$, and the learned observer gains $\hat Q_1(s)$, $\hat Q_2(s)$. The last row shows the error $Q_1(s)- \hat Q_1(s)$, $Q_2(s)- \hat Q_2(s)$. }
  		\label{fig:Q}
  	\end{figure*}

 Similar to the DeepONets for learning control kernels, except that the input channel for the first layer CNN is $2$.   Two DeepONets  are employed to learn the gain functions $Q_1$ and $Q_2$, respectively, each containing $6926913$ parameters.  
  	The two observer networks are trained together on  $1600$ instances, which only takes around 4 minutes. 
  	Figure~\ref{fig:Q} shows the analyzed observer gains, the learned DeepONet observer gains and the errors between them as $h=0.5s$.  
  	The network achieves a training loss of $5.02\rm{E}-5$ and a testing loss of $6.36\rm{E}-5$ after $300$ epochs, which are shown in Figure~\ref{fig:loss} (c).

\begin{figure*}   
  		\centering
  		\includegraphics[width=0.9\textwidth,trim=0 130 0 30,clip]{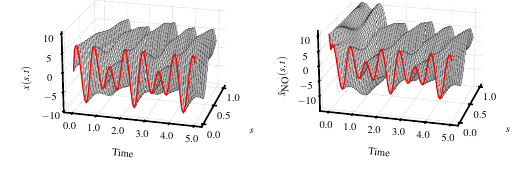} 
  		\caption{The system actual state $x(s,t)$ and the neural operator based observer $\hat x_{\mathrm{NO}}(s,t)$ under controller $U(t)=5\sin(3\pi t) +3\cos(2\pi t)$.  The left column showcases the actual state of the system state. 
   The right column showcases the estimated system state $\hat x_{\mathrm{NO}}(s,t)$ with the neural operator based observer.  Note  that the initial condition of the system is  $x(s,0)=\sin(2\pi s)$, while the initial condition of the neural operator observers is $\hat x_{\mathrm{NO}}(s,0)=10$.}
  		\label{fig:observer}
  	\end{figure*}

\begin{figure*}   
    \centering \includegraphics[width=0.9\textwidth,trim=0 130 0 30,clip]{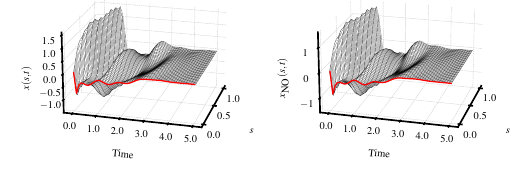}\\
    \begin{tabular}{cc}
\includegraphics[width=0.35\textwidth]{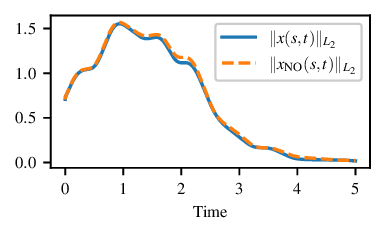}
&\includegraphics[width=0.4\textwidth]{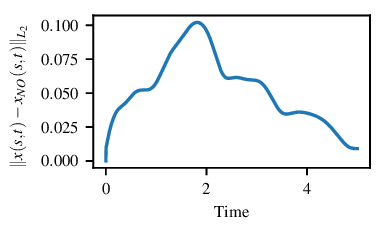}
    \end{tabular}
    \caption{ The closed-loop evolution under output feedback.  The left column in the fist row showcases the evolution of state $x(s,t)$ with the analytical kernels $K$, $L$, $J$ and observer gains $Q_1$, $Q_2$. 
The right column in the first row showcases the evolution of state $x_{\mathrm{NO}}(s,t)$ with NO kernels $\hat K$, $\hat L$, $\hat J$ and NO observer gains $\hat Q_1$, $\hat Q_2$.  The last row shows the $L_2$-norm of state $x$, $ x_{\mathrm{NO}}$, and the error between them.}
    \label{fig:output1}
\end{figure*}

  	Figure~\ref{fig:observer} demonstrates the convergence of the observation with the DeepONet learned gains to the system's actual state.  
  	In Figure~\ref{fig:output1}, we test the closed-loop system under the output feedback \eqref{estimated-controll} with three DeepONets approximating the control kernels and two DeepONets approximating the observer gains when the initial condition of the observer is set to the initial value of the system plus a random number that obeys $U(-1,1)$.  
  	
\begin{table*}[width=.75\textwidth]
  		\centering
  		\caption{Summary of kernel function calculation time consumption.}
  		\label{tab:time}       
  		\begin{tabular} {@{}ccccc@{}}%
   \toprule     
   &Model & \makecell[c]{Average Calculation Time \\(sec)\\(spatial step size $\Delta s =0.02$)}    & \makecell[c]{Average Calculation Time \\(sec)\\(spatial step size $\Delta s =0.01$)}  & \makecell[c]{Average Calculation Time\\ (sec)\\(spatial step size $\Delta s =0.005$)} \\ 
   \multirow{3}{*}{ \makecell[c]{Control kernels\\($K$, $L$, $J$)}} &Numerical solver & 0.025 & 0.654 & 1.692 \\
   &Neural operators & 0.011 & 0.0139 & 0.0265 \\	&Speedups& $2.3\times$ & $47.1\times$  &  $63.8\times$  \\
   \midrule 
   \multirow{3}{*}{\makecell[c]{Observer gains\\($Q_1$, $Q_2$)}}&Numerical solver  & 0.030 & 0.079 & 0.348 \\
   &Neural operators  & 0.006 & 0.006 & 0.007\\ &Speedups  &$5\times$ & $13.2\times$  &  $49.7\times$  \\
   \bottomrule 
  		\end{tabular}
\end{table*}
  	
  	Table \ref{tab:time} presents a comprison of the time consumption for kernel functions when solved numerically versus that generated by the trained DeepONets, respectively.   The term `average calculation time' refers to the mean time taken over 100 runs.  The duration required by numerical solvers grows significantly with the increase in discrete spatial step size, namely sampling precision. Conversely, the computation time for NOs shows only a marginal increase with larger spatial step sizes. However, the loss defined in \eqref{SL1Loss} still maintains on the order of $10^{-5}$.  It is worth noting that the approximation accuracy $\varepsilon$ defined in Theorem \ref{Th-Nkernel} and \ref{Th-NObs} corresponds to the Euclidean norm (or 2-norm) in vector space for the spatially discretized function (See. e.g. Theorem 2 in \cite{lu2019deeponet}). In the other words, the  accuracy, in terms of $\varepsilon$, is the square root of the loss \eqref{SL1Loss}. From the simulation results, we find that 
the accuracy of the control kernels is on the order of $10^{-3}$  and  that of the observer gains on the order of $10^{-2}$.

  	
  	\section{Conclusion}\label{6}

  	In this paper, we  apply the DeepONet operators to learn the  PDE backstepping  control kernels and observer of  a first-order hyperbolic PIDE system with state and sensor delays. Three neural operators are trained for the state feedback control from a group of numerical solutions of the backstepping kernel equations, which approximate three control kernel functions with the accuracy of magnitude of $10^{-3}$. The existence of arbitrary-precision NOs' approximation of the analytical kernel operators is proved by using the universal approximation theorem. The stability of the closed-loop system under state feedback with NO learning gains is also proved. Moreover, we use two neural operators to learn the observer gains and prove the observer with neural gains converge. The simulation results show that the accuracy of observer gains approximation can reach the magnitude of $10^{-2}$. Combined with the observer based  control system and observer error system, the stability of the output feedback system is proved, which verifies the separation principle under the neural operator gains. 
  	Further research will concern the delay-adaptive control of PDEs whose delays are unknown and control of high-dimensional PDEs control whose kernel functions are defined in higher spatial dimension.

\appendix
\section{Calculation of Kernel Function Equations}
Here is the calculation process to obtain \eqref{eq:K}--\eqref{eq:J}. We calculate the partial derivatives of transformation \eqref{eq:trans-ori} in time $t$ and space $s$ respectively, and we have
\begin{align}
z_t =&~ x_s+c(s) u(0,t)+\int_s^1 f(s,q) x(q,t) dq ~+K(s,1)x(1)\nonumber\\
   &-K(s,s)x(s)- \int_s^1 K_q(s,q) x(q)  dq \nonumber 
   \\&-\int_{s}^1 K(s,q)c(q) d q u(0,t )\\
   &~-\int^1_s K(s,q)
  	x(q,t)dq- \int_s^1 K(s,q)\int_q^1 f(q,r)x(r) dr dq \nonumber \\
   &~-L(s+h)v(1)+L(s) v(0) 
  	+h\int_0^1 L'(s+h r)v(r,t)dr\nonumber \\&~-J(s+\eta)u(1)+J(s)u(0)+\eta\int_0^1 J'(s+\eta r)u(r,t)dr,  \nonumber 
  \end{align}
  and
   \begin{align}
  	z_s(s,t) =&~ x_s(s,t)+K(s,s)x(s,t)-\int^1_s K_s(s,q)
  	x(q,t)dq\nonumber \\&~
  	-h\int_0^1 L'(s+h r )v(r,t)dr\nonumber \\&~-\eta\int_0^1 J'(s+\eta r)u(r,t)dr. 
  \end{align}
 Substitute them back into \eqref{eq:main-tar-z}--\eqref{eq:bnd-u-t}, then obtain the the kernel function equation \eqref{eq:K}--\eqref{eq:J}.   
   
	%
	%
	
	\bibliographystyle{cas-model2-names}
	
	\bibliography{Learningkernel}
\end{document}